\input amstex
\input xy
\xyoption{all}
\documentstyle{amsppt}
\document
\magnification=1200
\NoBlackBoxes
\nologo
\hoffset1.5cm
\voffset2cm
\vsize15.5cm

\def\G{\bold{G}}
\def\K{\bold{K}}
\def\N{\bold{N}}
\def\C{\bold{C}}
\def\Q{\bold{Q}}
\def\Z{\bold{Z}}
\def\R{\bold{R}}
\def\F{\bold{F}}

\def\A{\Cal{A}}
\def\S{\Cal{S}}
\def\V{\Cal{V}}
\def\H{\Cal{H}}

\def\Hom{\roman{Hom}\,}
\def\Spec{\roman{Spec}\,}
\def\Sets{\roman{Sets}\,}
\def\Burn{\roman{Burn}\,}
\def\Aut{\roman{Aut}\,}


\bigskip


\centerline{\bf  HOMOTOPY TYPES}

\medskip

\centerline{\bf  AND GEOMETRIES BELOW $\roman{Spec}\,\bold{Z}$}

\bigskip

\centerline{\bf Yuri I.~Manin, Matilde Marcolli}

\bigskip

{\bf Abstract.} After the first heuristic ideas about ``the field of one element'' $\bold{F}_1$
and ``geometry in characteristics 1'' (J.~Tits, C.~Deninger, M.~Kapranov, A.~Smirnov et al.),
there were developed several  general approaches to the construction
of ``geometries below $\roman{Spec}\,\bold{Z}$''. Homotopy theory and the ``the brave new algebra''
were taking more and more  important places  in these developments, systematically
explored by B.~To\"en and M.~Vaqui\'e, among others.

This article contains a brief survey and some new results  on {\it counting problems}
in this context, including various approaches to zeta--functions and generalised
scissors congruences.

\bigskip

 
  \medskip


\bigskip

\centerline{\bf  1. Brief summary and plan of exposition}

\medskip

{\bf 1.1. Geometries below  $\roman{Spec}\,\bold{Z}$: a general categorical framework.} Following
[ToVa09], Sec.~2.2 -- 2.5, we start with a symmetric monoidal category with unit 
$(C,\otimes , \bold{1})$. 
\smallskip
The category of commutative associative unital monoids $Comm\,(C)$ 
will play the role of commutative rings; accordingly, the opposite category $Aff_C := Comm\,(C)^{op}$
 will be an analogue of the category of {\it affine schemes.}
 \smallskip
 
 In order to be able to define more general schemes, objects of a category $Sch_C$, 
 we must introduce upon $Aff_C$
 {\it a Grothendieck topology} by giving a collection of sieves (covering families)
  defined for each object of  $Aff_C$.
  
  \smallskip
  
  It is shown in [ToVa09] that if  $(C,\otimes , \bold{1})$ is complete, cocomplete and closed,
  then there are several natural topologies upon $Aff_C$, whose names
  encode the similarities between them and respective topologies
  on the category of usual affine schemes (spectra of commutative rings),
  in particular, {\it Zariski topology.} Starting with monoidal symmetric category of
  abelian groups $(\bold{Z}$--Mod, $\otimes_{\bold{Z}}, \bold{Z})$ one comes to
  the usual category of schemes $Sch_{\bold{Z}}$ . 
  
  \smallskip
  
  The version of $\bold{F}_1$--schemes $Sch_{{\bold{F}}_1}$ suggested in [ToVa09]
  is embodied in the final stretch of the similar path starting with 
  $(Sets, \times , \{*\})$.

\smallskip

Finally, some pairs of categories $(Sch_C$, $Sch_D)$ can be related by two
functors going in reverse directions and satisfying certain adjointness properties.
Intuitively, one of them describes appropriate {\it descent data} upon certain objects $X_C$ 
of $Sch_C$ necessary and sufficient to define an object $X_D$ lying {\it under} $X_C$.
One of the most remarkable examples of such descent data from
$Sch_{{\bold{Z}}}$ to $Sch_{{\bold{F}}_1}$ was developed 
by J.~Borger, cf.  [Bo11a], [Bo11b]: roughly speaking, it consists in 
lifting all Frobenius morphisms to the respective $\bold{Z}$--scheme.
In a weakened form, when only a subset of Frobenius morphisms
is lifted, it leads to geometries  below $\roman{Spec}\,\bold{Z}$
but not necessarily over $\bold{F}_1$.

\medskip

{\bf 1.2. Schemes in the brave new algebra.} In a very broad sense,
the invasion of homotopy theory to mathematics in general started
with radical enrichment of Cantorian intuition about what are natural numbers $\N$:
they are cardinalities of not just arbitrary finite sets, but rather of sets of connected
components $\pi_0$ of topological spaces. 

\smallskip

The multiplication and addition in $\N$ have then natural lifts to
the world of stable homotopy theory, the ring of integers being enriched by passing to
sphere spectrum, where it becomes the initial object, in the same way
as $\Z$ itself is an initial object in the category of commutative rings, etc.
More details and references are given in the Section 4 of this article.

\smallskip

Moreover, ``counting functions'', such as numbers of $\F_q$--points of
a scheme reduced modulo a prime $p$, with $q=p^a$, can be generalized
to the world of scissors congruences where they become the basis for study of
zeta--functions.

\medskip

{\bf 1.3. The structure of the article.} 
(A). Sec.~2 is dedicated to a 
categorification in homotopy theory of a class of schemes in characteristics 1
admitting the following intuitive description:  {\it 1--Frobenius morphisms acting
upon their cohomology have eigenvalues that are roots of unity.}

\smallskip
Here are some more details.
The main arithmetic invariant of an algebraic
manifold $\V$ defined over a finite field $\F_q$, is its zeta--function $Z(\V,s)$
counting the numbers of its points $\roman{card}\,V(\F_{q^n})$
over finite extensions of $\F_q$.

\smallskip

Assuming for simplicity that $V$ is irreducible and smooth, 
we can identify $Z(V,s)$ with a rational function of $q^{-s}$
which is an alternating product of polynomials whose roots are
characteristic numbers of the Frobenius endomorphism $Fr_{q}$ of $V$
acting upon \'etale cohomology $H^*_{et}(V)$ of $V$.

\smallskip

In various versions of $\F_1$--geometry, the structure consisting
of cohomology with action of a Frobenius upon it is conspicuously
missing, although it is clearly lurking behind the scene
(see e.~g. a recent survey and study [LeBr17]).
\smallskip

Our homotopical approach here develops the analogy between Frobenius maps and Morse--Smale
diffeomorphisms.

\smallskip

In the remaining part of the article, we focus upon another "counting formalism" and its
well developed homotopical environment.

\smallskip
Namely, on the  $\F_q$-- schemes of finite type, non--necessarily smooth and proper  ones,  the function
$c:\, \V \mapsto \roman{card}\, \V(\F_q)$ satisfies the following ``scissors identity'':
if $X\subset \V$ is a Zariski closed subscheme, $Y:= \V\setminus X$, then 
$$
c(\V ) =c(X)+c(Y) 
$$
and also $c(\V_1\times \V_2)= c(\V_1) c(\V_2).$  In other words, $c$ becomes 
a ring homomorphism of the Grothendieck ring $K_0$ of the category of $\F_q$--schemes.

\medskip

(B). In Sec.~3 we consider a lift to the equivariant Grothendieck ring $K_0^{\hat\bold{Z}}(\Cal{V})$
of the integral Bost--Connes system, described in [CCMar09] in the context
of $\bold{F}_1$-geometry. We obtain endomorphisms $\sigma_n$ and additive
maps $\tilde\rho_n$ that act on $K_0^{\hat\bold{Z}}(\Cal{V})$ and map to the corresponding 
maps of the integral Bost--Connes system through the equivariant Euler characteristic.
We obtain in this way a noncommutative enrichment $\bold{K}_0^{\hat\bold{Z}}(\Cal{V})$ of the
Grothendieck ring and an Euler characteristic, which is a ring homomorphism
to the integral Bost--Connes algebra. After passing to $\bold{Q}$--coefficients, both algebras
become semigroup crossed products by the multiplicative semigroup of positive integers,
$$
 \bold{K}_0^{\hat\bold{Z}}(\Cal{V})\otimes_\bold{Z} \bold{Q} = K_0^{\hat\bold{Z}}(\Cal{V})_\bold{Q} \rtimes \bold{N} , 
 $$
with $K_0^{\hat\bold{Z}}(\Cal{V})_\bold{Q}=K_0^{\hat\bold{Z}}(\Cal{V})\otimes_\bold{Z} \bold{Q}$ and target of the
Euler characteristic the rational Bost--Connes algebra $\bold{Q}[\bold{Q}/\bold{Z}]\rtimes \bold{N}$.

\medskip

(C). In Sec.~4,  we  revisit the construction of the previous section, by further lifting it
from the level of the Grothendieck ring to the level of spectra. We use
the approach based on assemblers, developed in [Za17a, [Za17b].

\smallskip

We recall briefly the general formalism of [Za17a]  and  present
a small modification of the construction of [Za17c] of the assembler and
spectrum associated to the Grothendieck ring of varieties, which will be
useful in the following, namely the case of the equivariant Grothendieck
ring $K^{\hat\Z}(\V)$ considered in the previous section.

\smallskip

We then prove that the lift of the integral Bost--Connes algebra to the
level of the Grothendieck ring described in the previous section 
can further be lifted to the level of spectra. 

\medskip

(D). In Sec.~5 we discuss the construction of quantum statistical mechanical expectation values
on the Grothendieck ring $K_0^{\hat\Z}(\Cal{V})$ based on motivic
measures and the expectation values of the Bost--Connes system.

\medskip

(E). Finally, in Sec.~6  we revisit the construction discussed in the previous section in a setting where, 
instead of considering varieties with a good $\hat\Z$-action, we consider a ``dynamical"
model of $\F_1$-structure based on the existence of an endomorphism $f: X\to X$
that induces a quasi--unipotent map $f_*$ in homology. Our purpose here is to show
a compatibility between this proposal about $\F_1$-structures and the idea of [CCMar09]
of $\F_1$-geometry encoded in the structure of the integral Bost--Connes algebra,
through its relation to cyclotomic fields.

\smallskip

This also returns us to the framework of Sec.~2, thus closing the circle.

\bigskip

\centerline{\bf 2. Roots of unity as Weil numbers}

\medskip

{\bf 2.1. Local zetas and homotopy.}  In this section we are developing the idea,
sketched in the subsection 0.2 of [Ma10]: namely, that a $q=1$
replacement of the structure $(Fr_q, H^k_{et}(V))$
is a pair $(f_{*k}, H_k(M,\bold{Z}))$ where $f_{*k}$
is the action of a Morse--Smale diffeomorphism $f:\,M\to M$ 
of a compact manifold $M$ upon its homology $H_k(M,\bold{Z})$.

\smallskip

In particular, in this model characteristic roots of $f_*$ acting upon (co)homology
groups are roots of unity,  as might be expected from ``Weil numbers in
characteristic 1''.

\smallskip

We start with basic definitions and references.
\medskip

{\bf 2.2.~Morse--Smale diffeomorphisms.} Below $M$ always means a compact
smooth manifold, and $f:\,M\to M$ its diffeomorphism. A point $x\in M$ is called
non--wandering, if for any neighborhood $U$ of $x$, there is $n>0$ such that
$U\cap f^n(U)\ne \emptyset$.

\medskip

{\bf 2.2.1. Definition.} {\it $f$ is called a Morse--Smale diffeomorphism, if

(i) The number of non--wandering points of $f$ is finite.

(ii) $f$ is structurally stable that is, any small variation of $f$ 
is isotopic to $f$.}

\smallskip

This short definition appears in [Gr81]; for a more detailed discussion
of geometry, see [FrSh81].

\medskip

{\bf 2.3.~Action of Morse--Smale diffeomorphisms on homology.} For any
compact manifold $M$, its homology groups $H_k(M, \bold{Z})$
are finitely generated abelian groups. Any diffeomorphism
$f:\,M\to M$ induces automorphisms $f_{*k}$ of these groups.
According to [ShSu75], if $f$ is Morse--Smale, then each
$f_{*k}$ is quasi--unipotent, that is, its eigenvalues are roots of unity.

\smallskip

However, generally this condition is not sufficient. An additional necessary
condition is vanishing of a certain obstruction, described in [FrSh81]
and further studied in [Gr81].

\smallskip

Namely, consider the category  $QI$ whose objects are
pairs $(g,H)$, where $H$ is a (finitely generated) abelian group,
and $g:\,H\to H$ is a quasi--idempotent endomorphism of $H$
(by definition, this means that eigenvalues of $g$ are zero
or roots of unity.)  Morphisms in $QI$ are self--evident.

\smallskip

The abelian group $K_0(QI)$ admits a morphism onto its torsion subgroup $G$.
 Denote by $\varphi:\, K_0(QI)\to G$ be such a morphism.

\medskip

{\bf 2.3.1. Theorem ([FrSh81]).} {\it Let $f:\,M\to M$ be a diffeomorphism such that
each $f_{*k}:\, H_k(M,\bold{Z}) \to H_k(M,\bold{Z})$ is quasi--unipotent.
Let $[f_{*k}]$ be the class in $K_0(QI)$ of the pair $(f_{*k},H_k(M,\bold{Z})).$
 Then
\smallskip
(i) If $f$ is Morse--Smale, then  $\chi(f_*):= \sum_k (-1)^k \varphi ([f_{*k}])\in G$
is zero.
\smallskip

(ii) If in addition $M$ is simply connected, and $\roman{dim}\,M > 5$,
then vanishing of $\chi(f_*)$ implies that $f$ is isotopic to a Morse--Smale
diffeomorphism.}

\medskip

{\bf 2.4.~ Passage from $(f_{*k}, H_k(M,\bold{Z}))$ to a scheme over $F_1.$}
In this subsection, we sketch a final step from homotopy to $F_1$--geometry.

\smallskip

One should keep in mind, however, that there might be many divergent paths,
starting at this point, because there are several different versions of 
``geometries over $F_1$''.

\smallskip

We will choose here a version developed in James Borger's paper
[Bo11a] (and continued in [Bo11b]). Roughly speaking, in order to define an affine
scheme over $F_1$, one should give a (commutative) ring with $\lambda$--structure
and then treat this $\lambda$--structure as ``descent data'' from the base $\roman{Spec}\,\bold{Z}$
to the base $\roman{Spec}\,F_1.$

\medskip

{\bf 2.4.1.~ From $(f_{*k}, H_k(M,\bold{Z}))$ to a $\lambda$--ring.} Consider for simplicity
only the case (ii) of Theorem 3.1. Then for each $k$,  $H_k(M,\bold{Z}))$
is a free $\bold{Z}$--module of finite rank, and $f_{*k}$ is its quasi--unipotent automorphism.
Fix such a $k$.

\smallskip

Introduce upon  $H_k(M,\bold{Z}))$ the structure of $\bold{Z} [T,T^{-1}]$--module,
where $T$ acts as $f_{*k}$. 

\smallskip

We can consider the minimal subcategory $\Cal{C}$ of  $\bold{Z} [T,T^{-1}]$--modules,
containing $H_k(M,\bold{Z})$ and closed with respect to direct sums, tensor products,
and exterior powers, and then produce the Grothendieck $\lambda$--ring $K_0(\Cal{C})$
using exterior powers to define the relevant $\lambda$--structure (cf.~[At61], p.~26, and 
[Le81]).

\medskip

{\bf 2.4.2. Definition.} {\it The $F_1$--scheme, whose Borger's lift to $\roman{Spec}\,\bold{Z}$
is $K_0(\Cal{C})$, is the representative of  $(f_{*k}, H_k(M,\bold{Z}))$ in $F_1$--geometry.
}
\medskip

{\bf 2.4.3. Remark.} It seems that another short path from  $(f_{*k}, H_k(M,\bold{Z}))$ 
to an $F_1$--scheme defined differently might lead to one of Le Bruyn's  spaces
in [LeBr17].

\medskip

{\bf 2.5. Cases when eigenvalues of conjectural Frobenius maps are not
roots of unity.} Here I will briefly discuss possible extensions of the picture above,
leading to various geometries  ``below $\roman{Spec}\,\bold{Z}$'' but generally
{\it not} over $\roman{Spec}\,F_1$.

\smallskip

The most interesting new virtual zeta--functions of this type
were discovered only recently under the generic name ``zeta--polynomials'':
cf.~[Ma16], [JiMaOnSo16], [OnRoSp16].

\smallskip

In [Ma16], it was described how to produce such polynomials from period
polynomials of any cusp form over $SGL(2,\bold{Z})$ which is an eigenform for
all Hecke operators: this passage is a kind of ``discrete Mellin transform''. 
It was also proved that zeroes of period polynomials lie
on the unit circle of the complex plane, but generally are {\it not}
roots of unity.

\smallskip

Both this construction and the results about zeroes were generalised  in 
[JiMaOnSo16], [OnRoSp16] to the case of  cusp newforms of even weight for $\Gamma_0(N)$.
It turned out that, with appropriate scaling, zeroes of period polynomials lie on the circle
$\{z\,|\, |z|=\dfrac{1}{\sqrt{N}}\}$.

\smallskip

{\bf 2.5.1. Problem.} Make explicit geometries under $\roman{Spec}\,\bold{Z}$ in which one
can accommodate 
the respective
``motives''.

\medskip

\centerline{\bf 3. The Bost--Connes system and the Grothendieck ring}

\medskip

{\bf 3.1. $\hat\bold{Z}$-equivariant Grothendieck ring ([Lo99]).}
We recall the following definitions from [Lo99]. Let  $G$ be an affine algebraic group
acting upon an algebraic  variety $X$. We say that this action
is {\it good}, if
each $G$--orbit is contained in an affine open subset of $X$.  
\smallskip

The Grothendieck ring
$K_0^G(\Cal{V})$ is generated by isomorphism classes $[X]$ of pairs, consisting of varieties with good $G$--action.
Upon these pairs
 the inclusion--exclusion relations are imposed: $[X]=[Y]+[X\smallsetminus Y]$ 
where $Y\hookrightarrow X$ and $X\smallsetminus Y \hookrightarrow X$
are $G$--equivariant embeddings. Multiplication in $K_0^G(\Cal{V})$ 
is induced by the diagonal $G$--action on the product.

\smallskip

In the main special case considered in [Lo99], 
$$ 
G=\hat\bold{Z}=\varprojlim_n \bold{Z}/n\bold{Z} .
$$ 
One defines $K_0^{\hat\bold{Z}}(\Cal{V})$
as the Grothendieck ring of varieties with an action of $\hat\bold{Z}$
that factors through a good action of some finite quotient $\bold{Z}/n\bold{Z}$.

\smallskip

In this section, we first consider varieties defined over the field $\C$ of complex numbers 
and classes in the Grothendieck ring correspondingly taken in $K_0(\V_\C)$
and the equivariant $K_0^{\hat\Z}(\V_\C)$. We then consider the case of varieties
defined over $\Q$ with the equivariant Grothendieck ring $K_0^{\hat\Z}(\V_\Q)$.
In the first case the target of the equivariant Euler characteristic consists of the abelian
part $\Z[\Q/\Z]$ of the integral Bost--Connes algebra, while in the second case
it is a subring spanned by the range projectors of the Bost--Connes algebra.

\smallskip

As observed in [Lo99], there is an Euler characteristic ring homomorphism
$$ 
\chi^{\hat\Z}: K_0^{\hat\Z}(\V_\C) \to K_0^{\hat\Z}(\C), 
$$
to the Grothendieck ring of finite dimensional representations. Since the
character group is $\Hom(\hat\Z,\G_m)=\Q/\Z$, the latter is identified
with the group ring
$$ 
K_0^{\hat\Z}(\C) = \Z [\Q/\Z]. 
$$

\medskip

{\bf 3.2. Equivariant Euler characteristics.} Since we are considering varieties with a good action of $\hat{\Z}$
that factors through some finite quotient, the Euler characteristic above can be
obtained as in [Ve73], for an action of a finite group $G$, by replacing 
the alternating sum of dimensions of the cohomology groups with a sum in 
the ring $R(G)$ of representations of $G$ of the classes of cohomology groups, 
viewed as $G$-modules.

\smallskip

As observed in [Gu-Za17], one can also define an equivariant
Euler characteristic, for good actions on varieties of a finite group $G$, 
as a ring homomorphism $\chi^G: K_0^G(\V) \to A(G)$, where
$A(G)$ is the Burnside ring of $G$, the Grothendieck ring of the
category of finite $G$-sets. In this case the equivariant Euler characteristic
is defined as $\chi^G(X)=\sum_{k\geq 0} [X_k]$ where 
$X$ has a simplicial decomposition with $k$--skeleton $X_k$, 
such that $G$ acts by simplicial maps
which map each $k$--simplex either identically to itself or to another simplex,
so that it makes sense to consider the classes $[X_k]$ in $A(G)$. It is
shown in [Gu-Za17] that the result is independent of the choice of such a
simplicial decomposition. It is also shown that any invariant with values in
a commutative ring $R$, defined on varieties with a good $G$--action
homeomorphic to locally closed unions of cells in finite CW-complexes 
with $G$ acting by cell maps, which satisfies inclusion-exclusion (on
$G$-invariant decompositions) and multiplicativity on products is
necessarily a composition 
$$ 
K_0^G(\V) \to A(G) \to R 
$$
of $\chi^G$ with a ring homomorphism $\varphi: A(G)\to R$.
In particular, this is the case for the Euler characteristic $K_0^G(\V) \to R(G)$
obtained by composing $\chi^G: K_0^G(\V) \to A(G)$ with 
the natural ring homomorphism $\varphi: A(G) \to R(G)$ that sends $G$--sets
to their space of functions. 

\smallskip

When considering the profinite group $\hat\Z =\varprojlim_n \Z/n\Z$, the
Burnside rings $A(\Z/n\Z)$ form a projective system with limit
$$ \hat A(\hat\Z) = \varprojlim_n A(\Z/n\Z), $$
the completed Burnside ring of $\hat\Z$, 
which is the Grothendieck ring of almost finite $\hat\Z$-spaces,
namely those $\hat\Z$-spaces that are discrete and essentially
finite, that is, such that for every open subgroup $H$ the set
of points fixed by all elements of $H$ is finite, see Section 2 of [DrSi88]. 
The Burnside ring $A(\hat\Z)$ of finite $\hat\Z$-spaces sits as a dense subring
of $\hat A(\hat\Z)$. Moreover, there is an identification of this completed Burnside ring
with the Witt ring $\hat A(\hat\Z)= W(\Z)$, see Corollary~1 of [DrSi88].

\medskip
{\bf 3.3. Lifting the integral Bost--Connes system.}
Consider now the endomorphisms $\sigma_n: \Z [\Q/\Z] \to \Z[\Q/\Z]$ given
by $\sigma_n(e(r))=e(nr)$ on the standard basis $\{ e(r)\,:\, r\in \Q/\Z \}$ of
$\Z [\Q/\Z]$. The integral Bost--Connes algebra $\A_\Z$ introduced in [CCMar09]
is generated by the group ring $\Z[\Q/\Z]$ together with elements $\tilde\mu_n$
and $\mu_n^*$ satisfying the relations
$$
\tilde\mu_{nm}=\tilde\mu_n \tilde\mu_m, \ \ \ \ \mu_{nm}^*=\mu_n^* \mu_m^*, \ \ \ \ 
\mu_n^* \tilde \mu_n =n, \ \ \ \ \tilde\mu_n \mu_n^* =\mu_n^* \tilde\mu_n, 
\eqno(3.1)
$$
 where the first two relations hold for arbitrary $n,m\in \N$, the third for arbitrary $n\in \N$ and
the fourth for $n,m\in \N$ satisfying $ (n,m)=1$, and the relations
$$
x \tilde\mu_n = \tilde\mu_n \sigma_n(x) \ \ \ \ \mu_n^* x = \sigma_n(x) \mu_n^*, \ \ \ \ 
\tilde\mu_n  x \mu_n^* = \tilde\rho_n(x),
\eqno(3.2)
$$
for any $x\in \Z[\Q/\Z]$, where $\tilde\rho_n(e(r))=\sum_{nr'=r} e(r')$.

\smallskip

The maps $\tilde\rho_n$ and the endomorphisms $\sigma_n$ satisfy the
compatibility conditions, for all $x,y\in \Z[\Q/\Z]$ (see Proposition~4.4 of [CCMar09])
$$
\tilde\rho_n(\sigma_n(x)y) = x \tilde\rho_n(y), \ \ \ \ \  \sigma_n(\tilde\rho_m(x))=
(n,m) \cdot \tilde\rho_{m'}(\sigma_{n'}(x)),
\eqno(3.3)
$$
where $(n,m)=\gcd\{n,m\}$ and $n'=n/(n,m)$ and $m'=m/(n,m)$.

\medskip
{\bf 3.3.1. Lemma.} {\it
The endomorphisms $\sigma_n: \Z [\Q/\Z] \to \Z[\Q/\Z]$ lift to endomorphisms
$\sigma_n : K_0^{\hat\Z}(\V_\C) \to K_0^{\hat\Z}(\V_\C)$ such that the
following diagram commutes}
$$ 
\xymatrix{  K_0^{\hat\Z}(\V_\C) \ar[d]^{\sigma_n} \ar[r]^{\chi^{\hat\Z}} & \Z [\Q/\Z] \ar[d]^{\sigma_n} \\ K_0^{\hat\Z}(\V_\C) \ar[r]^{\chi^{\hat\Z}} & \Z [\Q/\Z].}
 $$
{\it These endomorphisms define a semigroup action of the multiplicative semigroup $\N$ on the
Grothendieck ring $K_0^{\hat\Z}(\V_\C)$. }

\medskip

{\bf Proof.} Let $X$ be a variety with a good $\hat\Z$-action, which factors through some
finite quotient $\Z/N\Z$. Let $\alpha: \hat\Z \times X \to X$ denote the action. The
endomorphisms $\sigma_n: \Z[\Q/\Z] \to \Z[\Q/\Z]$ given by $\sigma_n(e(r))=e(nr)$
have an equivalent description as the action on the group of roots of unity of all orders given by raising to 
the $n$--th power $\sigma_n: \zeta \mapsto \zeta^n$.
One can then obtain an action $\alpha_n: \hat\Z \times X \to X$ given by $\alpha_n = \alpha\circ \sigma_n$.
Thus, we assign to a pair $(X,\alpha)$ of a variety with a good $\hat\Z$--action the pair $(X,\alpha_n)$
of the same variety with the action $\alpha_n$. This assignment respects isomorphism classes and is
compatible with the relations, hence it determines endomorphisms of $K_0^{\hat\Z}(\V_\C)$,
with $\sigma_{nm}=\sigma_n\circ \sigma_m$, namely a semigroup action of the multiplicative semigroup $\N$
on the Grothendieck ring $K_0^{\hat\Z}(\V_\C)$. $\blacksquare$

\medskip

The maps $\tilde\rho_n: \Z[\Q/\Z]\to \Z[\Q/\Z]$
of the form 
$$
\tilde\rho_n(e(r))=\sum_{nr'=r} e(r') 
\eqno(3.4)
$$
are not ring homomorphisms but only morphisms
of abelian groups. After tensoring with $\Q$, one obtains the group algebra $\Q[\Q/\Z]=\Z[\Q/\Z]\otimes_\Z \Q$
on which the $\tilde\rho_n$ induce endomorphisms of the form
$\rho_n(e(r))=n^{-1} \sum_{nr'=r} e(r')$ satisfying the relations $\sigma_n \rho_n(x)=x$ and
$\rho_n\sigma_n(x)=\pi_n x$, for $x\in \Q[\Q/\Z]$ and the idempotent $\pi_n=n^{-1}\sum_{ns=0} e(s)$.
The arithmetic Bost--Connes algebra is the crossed product $\Q[\Q/\Z]\rtimes \N$ generated by
$\Q[\Q/\Z]$ and $\mu_n$, $\mu_n^*$ with the crossed-product action of $\N$ implemented 
by $\mu_n x \mu_n^*=\rho_n(x)$, see [CCMar09]

\smallskip

Once one considers varieties defined over $\Q$, the Grothendieck ring $K_0^{\hat\Z}(\Q)$ 
can be characterized as follows.
\medskip
{\bf 3.3.2. Lemma.} {\it The Grothendieck ring $K_0^{\hat\Z}(\Q)$ can be identified with the
subring of $\Z [\Q/\Z]$ generated by the elements $n\pi_n = \sum_{ns=0} e(s)$.}

\medskip

{\bf Proof.} As noted in [Lo99], the 
element $\sum_{ns=0} e(s)$ in $\Z[\Q/\Z]$ is the image of the irreducible representation 
of $\Z/n\Z$ given by the cyclotomic field $\Q(\zeta_n)$ seen as a $\Q$-vector space, and
these representations give a basis of $K_0^{\hat\Z}(\Q)$.  $\blacksquare$

\medskip

{\bf 3.3.3. Remark.} According to the Corollary~4.5 in [CCMar09], the range of the 
maps $\tilde\rho_n$ in (3.4) is an ideal in $\Z [\Q/\Z]$. This
follows from the relation $\tilde\rho_n(\sigma_n(x)y)=x \tilde\rho_n(y)$ of (3.3).
If $r'$ an element of the set 
$E_n(r)=\{ r'\in\Q/\Z\,:\, nr' =r \}$, we have 
$\tilde\rho_n(e(r))=e(r') \sum_{ns=0} e(s)$.

\medskip

{\bf 3.3.4. Lemma.} {\it There are endomorphisms $\sigma_n: K_0^{\hat\Z}(\Q) \to K_0^{\hat\Z}(\Q)$
induced by the endomorphisms $\sigma_n: \Z[\Q/\Z]\to \Z[\Q/\Z]$.
They lift to endomorphisms $\sigma_n : K_0^{\hat\Z}(\V_\Q) \to 
K_0^{\hat\Z}(\V_\Q)$ given by $\sigma_n(X,\alpha)=(X,\alpha\circ \sigma_n)$
as in Lemma 3.3.1.
}

\medskip

{\bf Proof.} By identifying $K_0^{\hat\Z}(\Q)$ with a subring of $\Z[\Q/\Z]$ 
as in Lemma 3.3.2.
 we see that the endomorphisms 
$\sigma_n:  \Z[\Q/\Z]\to \Z[\Q/\Z]$ induce endomorphisms of $K_0^{\hat\Z}(\Q)$
by the relations $\sigma_n (\tilde\rho_m(x))=(n, m) \cdot \tilde\rho_{m'}(\sigma_{n'}(x))$
as in (3.3.3). Since we have $n\pi_n =\tilde\rho_n(1)$,
we see that the endomorphisms $\sigma_n$ map the subring $K_0^{\hat\Z}(\Q)$
to itself. The lift to $\sigma_n : K_0^{\hat\Z}(\V_\Q) \to K_0^{\hat\Z}(\V_\Q)$ is obtained by 
by the same argument as in Lemma 3.3.1. Namely, the map 
$\sigma_n(X,\alpha)=(X,\alpha\circ \sigma_n)$ defines an endomorphism of 
$K_0^{\hat\Z}(\V_\Q)$ which satisfies $\sigma_n \circ \chi^{\hat\Z}=\chi^{\hat\Z} \circ \sigma_n$,
with $\chi^{\hat\Z}:  K_0^{\hat\Z}(\V_\Q) \to K_0^{\hat\Z}(\Q)$ the Euler characteristic.
$\blacksquare$

\medskip

In order to lift the maps $\tilde\rho_n$ to the level of the Grothendieck
ring $K_0^{\hat\Z}(\V_\Q)$, we put
$$
\bold{V}_n := [Z_n,\gamma_n] \,\, \in K_0^{\hat\Z}(\V_\Q), \ \ \  
\eqno(3.5)
$$
where $Z_n$ is a zero--dimensional variety over $\Q$ with $\# Z_n(\bar \Q)=n$. Any such
variety with a smooth model over $\Spec \Z$ can be identified with $Z_n=\Spec(\Q^n)$. It
is endowed with the natural action $\gamma_n: \hat\Z \times Z_n \to Z_n$ that factors 
through $\Z/n\Z$.

\smallskip

Given a variety $X$ with a good $\hat\Z$-action $\alpha: \hat\Z \times X \to X$,
let $$
\Phi_n(\alpha): \hat\Z \times X \times Z_n \to X \times Z_n
$$ 
be 
given by
$$
\Phi_n(\alpha)(\zeta, x, a_i) =  (x,\gamma_n(\zeta,a_i))\quad \roman{for}\quad i=1,\ldots,
$$
$$
\roman{and}\quad
(\alpha(\zeta,x),\gamma_n(\zeta,a_n))\quad  \roman{for}\quad i=n.
\eqno(3.6)
$$

Notice that (3.6) is just a form of the {\it Verschiebung} map: for $\zeta$ 
the generator of $\Z/n\Z$ we have

$$
\Phi_n(\alpha)(\zeta, x, a_i) =  (x,a_{i+1})\quad \roman{for}\quad i=1,\ldots, n-1, \quad
$$
$$
\roman{and}\quad
(\alpha(\zeta,x), a_1)\quad  \roman{for}\quad i=n.
$$
\medskip

{\bf 3.3.5. Proposition.} {\it The maps
$$
\tilde\rho_n [X,\alpha] := [X\times Z_n, \Phi_n(\alpha)]
\eqno(3.7)
$$
define a homomorphism of the Grothendieck group $K_0^{\hat{\Z}}(\Cal{V}_{\Q})$
that satisfies
$$
\sigma_n\circ \tilde\rho_n [X,\alpha] = [X,\alpha]^{\oplus n}
\eqno(3.8)
$$
and
$$
 \tilde\rho_n \circ \sigma_n [X,\alpha] = \tilde\rho_n [X,\alpha\circ\sigma_n] =
 [X,\alpha]\cdot [Z_n,\gamma_n].
\eqno(3.9)
$$
}

\smallskip

{\bf Proof.} Given a variety $X$ with a good $\hat\Z$-action $\alpha:\hat\Z\times X \to X$,
consider the product $[X,\alpha]\cdot \bold{V}_n$ in $K_0^{\hat\Z}(\V_\Q)$. This
class has a representative $[X\times Z_n, (\alpha\times \gamma_n)\circ \Delta]$,
where $\Delta: \hat\Z \to \hat\Z\times \hat\Z$ is the diagonal.
We have
$$ \sigma_n \circ \tilde\rho_n [X,\alpha]=\sigma_n [X\times Z_n, \Phi_n(\alpha)]=
[X\times Z_n, \Phi_n(\alpha)\circ \sigma_n]. $$
We also have
$$ \Phi_n(\alpha)\circ \sigma_n = (\alpha\times 1)\circ \Delta, $$
since 
$$ \Phi_n(\alpha)\circ \sigma_n (\zeta,x,a_i)=\Phi_n(\alpha) (\zeta^n,x,a_i) =
(\Phi_n(\alpha)(\zeta))^n (x,a_i), $$
where we write $\Phi_n(\alpha)(\zeta): X\times Z_n \to X\times Z_n$ for
the action of $\zeta\in \hat\Z$, with $\Phi_n(\alpha)(\zeta_1\cdots \zeta_n)=
\Phi_n(\alpha)(\zeta_1)\circ \cdots \circ \Phi_n(\alpha)(\zeta_n)$, and the $n$-fold
composition gives
$$ \Phi_n(\alpha)(\zeta) \circ \cdots \circ \Phi_n(\alpha)(\zeta)  (x,a_i)=(\alpha(\zeta,x),a_i). $$
This shows (3.8). 
\smallskip

The second relation is obtained similarly. We have
$$
 \tilde\rho_n \circ \sigma_n [X,\alpha] = \tilde\rho_n [X,\alpha\circ\sigma_n] =
[X\times Z_n, \Phi_n(\alpha\circ \sigma_n)], 
$$
where 
$$
\Phi_n(\alpha)(\zeta, x, a_i) =  (x,a_{i+1})\quad \roman{for}\quad i=1,\ldots, n-1, \quad
$$
$$
\roman{and}\quad
(\alpha(\zeta^n,x), a_1)\quad  \roman{for}\quad i=n.
$$
Now $\alpha(\zeta^n,x)=\alpha(\zeta)^n (x)$, hence we have
$\Phi_n(\alpha\circ \sigma_n)(\zeta) (x,a_i) = \Phi_n(\alpha(\zeta)^n) (x,a_i)$.
The usual relations $V_n(F_n(a)b)=a V_n(b)$ between Frobenius $F_n$ and
Verschebung $V_n$ (see Proposition~2.2 of [CC14]) holds in this case
as well in the form $\Phi_n(\alpha(\zeta)^n)=\alpha(\zeta) \Phi_n(1)$ where
$\Phi_n(1)(x,a_i)=(x,a_{i+1})$ is the cyclic permutation action of $\Z/n\Z$ on $Z_n$.
Thus, we obtain $\Phi_n(\alpha\circ \sigma_n) =(\alpha \times \gamma_n)\circ \Delta$.
This gives (3.9). $\blacksquare$

\medskip

The relation (3.8) corresponds to $\sigma_n \circ \tilde\rho_n(x)= n x$, 
and (3.9) to $\tilde\rho_n\circ \sigma_n (x) = n \pi_n x$, for $x\in \Z[\Q/\Z]$.
They are geometric manifestations of the same relation between the maps $\sigma_n$
and $\tilde\rho_n$ of the integral Bost--Connes system and the Frobenius and Verschiebung
described in [CC14]. For other occurrences of the same relation see also [MaRe17], [MaTa17].

\smallskip

{\bf 3.4.  A non--commutative extension of the Grothendieck ring.}  
Let $\K_0^{\hat\Z}(\V_\Q)$ be the non-commutative ring generated by 
$K_0^{\hat\Z}(\V_\Q)$ and elements $\tilde\mu_n$, $\mu_n^*$ for $n\in \N$
satisfying the relations (3.1) for all $n,m\in \N$, and (3.2)
for all $x=[X,\alpha] \in K_0^{\hat\Z}(\V_\Q)$ and all $n\in \N$.

\medskip

{\bf 3.4.1. Lemma.} {\it
The Euler characteristic $\chi^{\hat\Z}: K_0^{\hat\Z}(\V_\Q) \to K_0^{\hat\Z}(\Q) \hookrightarrow \Z[\Q/\Z]$
extends to a ring homomorphism $\chi: \K_0^{\hat\Z}(\V_\Q)\to \A_\Z$ to the integral Bost--Connes algebra.
After tensoring with $\Q$, we obtain a homomorphism of semigroup crossed product rings
$$ \chi: K_0^{\hat\Z}(\V_\Q)_\Q \rtimes \N \to \Q[\Q/\Z]\rtimes \N, $$
where $\K_0^{\hat\Z}(\V_\Q)\otimes_\Z \Q =K_0^{\hat\Z}(\V_\Q)_\Q \rtimes \N$ with
$K_0^{\hat\Z}(\V_\Q)_\Q= K_0^{\hat\Z}(\V_\Q)\otimes_\Z \Q$,
and $\A_\Z\otimes_\Z \Q =\A_\Q =\Q[\Q/\Z]\rtimes \N$ is the rational Bost--Connes algebra.
}

\medskip

{\bf Proof.}  We define the map $\chi$ as $\chi^{\hat\Z}$ on elements of $K_0^{\hat\Z}(\V_\Q)$
and the identity on the extra generators $\chi(\tilde\mu_n)=\tilde\mu_n$ and $\chi(\mu_n^*)=\mu_n^*$.
By Lemma 3.3.4 and Proposition 3.3.5, this map is compatible with the
relations in $\K_0^{\hat\Z}(\V_\Q)$ and in $\A_\Z$. After tensoring with $\Q$, the algebra
$\K_0^{\hat\Z}(\V_\Q)\otimes_\Z \Q$ can be identified with a semigroup crossed product by
taking as generators the elements of $K_0^{\hat\Z}(\V_\Q)$ and $\mu_n=n^{-1} \tilde \mu_n$
and $\mu_n^*$, which satisfy the relations
$$
 \mu_n^*\mu_n =1, \ \ \ \ \mu_{nm}=\mu_n\mu_m, \ \ \ \ \mu_{nm}^*=\mu_n^*\mu_m^*, \ \  \forall n,m\in\N, 
\ \ \mu_n \mu_m^* = \mu_m^*\mu_n \text{ if } (n,m)=1 $$
$$ \mu_n x \mu_n^* = \rho_n(x) \ \ \text{ with } \ \  \rho_n(x) =\frac{1}{n} \tilde\rho_n (x), 
$$
with $\sigma_n \rho_n(x) =x$, for all $x=[X,\alpha] \in K_0^{\hat\Z}(\V_\Q)$. 
The semigroup action in the crossed product $K_0^{\hat\Z}(\V_\Q)_\Q \rtimes \N$ is given by
$x \mapsto \rho_n(x) =\mu_n x \mu_n^*$. 
The target algebra is the rational Bost--Connes algebra $\A_\Z\otimes_\Z \Q =\Q[\Q/\Z]\rtimes \N$.
Again the map given by $\chi^{\hat\Z}$ on elements of $K_0^{\hat\Z}(\V_\Q)$ and by
$\chi(\mu_n)=\mu_n$ and $\chi(\mu_n^*)=\mu_n^*$ determines a homomorphism of crossed-product
algebras.  $\blacksquare$

\bigskip

\centerline{\bf  4. From Rings to Spectra}

\medskip

{\bf 4.1. Spectra.} We give a brief review of spectra, with the purpose of
recalling a construction of Segal [Se74] that constructs spectra from $\Gamma$--spaces. 
We then review the notion of assembler categories [Za17a] and how 
they can be used to construct a $\Gamma$--space and an associated spectrum
whose $\pi_0$ realises certain abstract scissor--congruence relations.

\smallskip

The construction of spectra from $\Gamma$--spaces was first developed
in the context of the Bousfield--Friedlander spectra, see Definition~2.1 of [BousFr78].

\smallskip

In this setting, one considers the simplicial category $\Delta$, which has an object $[n]$ 
for each $n\in \N$ given by the finite totally ordered set $[n]=\{ 0 < 1 < \ldots < n-1 \}$,
with morphisms the face and degeneracy maps $\delta^n_i$ and $\sigma^n_i$
satisfying the simplicial relations. 

\smallskip

A simplicial object is a contravariant functor $S:\Delta^{op}\to \Cal{C}$ from $\Delta$ to a
given category $\Cal{C}$. It is determined by a sequence of objects $X(n)$ of $\Cal{C}$ with
morphisms corresponding to faces and degeneracies. We denote by $\Delta(\Cal{C})$
the resulting category of simplicial objects in $\Cal{C}$. In particular, a simplicial set is a 
simplicial object in the category of sets and we will use the notation $\underline{\Delta}=\Delta({\Sets})$
for the category of simplicial sets. 

\smallskip
Similarly, a bisimplicial object is a functor $BS : \Delta^{op} \times \Delta^{op} \to \Cal{C}$,
or equivalently a simplicial object in the category of simplicial objects $\Delta(\Cal{C})$.
The diagonal of a bisimplicial object $BS$ is the simplicial object obtained by
precomposition  of $BS$ with the diagonal functor $\Delta^{op}\to\Delta^{op}\times \Delta^{op}$.

\smallskip

The category $\S$ of Bousfield--Friedlander spectra has objects $X$ given by
sequences of simplicial sets $X=\{ X_n \}_{n\geq 0}$ endowed with structure
maps $\varphi^X_n: S^1\land X_n \to X_{n+1}$ for all $n \geq 0$, and morphisms
given by maps $f_n: X_n \to Y_n$ with commutative diagrams
$$
 \xymatrix{ S^1\land X_n \ar[r]^{\varphi^X_n} \ar[d]^{1_{S^1}\land f_n} & X_{n+1}\ar[d]^{f_{n+1}} \\ 
S^1\land Y_n \ar[r]^{\varphi^Y_n} & Y_{n+1}. } 
$$

\smallskip

The sphere spectrum $\bold{S}$ has $\bold{S}_n =S^1\land \cdots \land S^1$, the $n$--fold smash product, 
and $\varphi_n$ the identity map. 

\smallskip

Let $\gamma^X_n: X_n \to \Omega X_{n+1}$ be the maps induced by the adjoints of 
the structure maps. An $\Omega$-spectrum is a spectrum where the maps $\gamma^X_n$ are
weak equivalences for all $n$.

\smallskip

The homotopy groups $\pi_k(X)$ of spectra are given by
$$ \pi_k(X) =\varinjlim_n \pi_{k+n} (X_n) $$
over the maps $\pi_{k+n} (X_n)  \to \pi_{k+n} (\Omega X_{n+1})\simeq \pi_{n+k+1} (X_{n+1})$,
induced the $\gamma^X_n$. 
A spectrum is $n$--connected if $\pi_k(X)=0$ for all $k\leq n$ and connective if it is $-1$-connected. 
A spectrum $X$ is cofibrant if all the structure maps $\phi^X_n: S^1\land X_n \to X_{n+1}$
are cofibrations. 

\smallskip

However, a problem with the Bousfield--Friedlander spectra is that they do not have 
a homotopically good smash product.  Constructions of categories of spectra with
smash products were developed in the '90s, in particular, the $S$--modules
model of [EKMM97] and the symmetric spectra model of [HSS00]. 
In a more modern approach, it is therefore preferable to work with 
symmetric spectra for the Segal construction. Indeed, the $\Gamma$--spaces, 
SW-categories and Waldhausen categories that occur in relation to
the spectra underlying the Grothendieck ring of varieties and its
variants naturally give rise to symmetric spectra.

\smallskip

A symmetric spectrum consists of a sequence of pointed spaces (pointed simplicial sets)
$X=\{ X_n \}_{n\geq 0}$ together with a left action of the symmetric group $S_n$ on
$X_n$ for all $n\geq 0$ and structure maps given by based maps 
$\varphi^X_n: S^1\land X_n \to X_{n+1}$ for all $n \geq 0$, with the condition that, for
all $n,m\geq 0$ the composition $\varphi^X_{n+m-1}\circ\cdots\circ \varphi^X_n$
$$ \varphi^X_{n+m-1}\circ\cdots\circ \varphi^X_n: S^m \land X_n \to S^{m-1}\land X_{n+1} \to \cdots S^1 \land X_{n+m-1}
\to X_{n+m} $$
is $S_n \times S_m$-equivariant.
\smallskip
A morphism of symmetric spectra is a collection of $S_n$--equivariant
based maps $f_n : X_n \to Y_n$ such that $f_{n+1}\circ \varphi^X_n = \varphi^Y_n \circ (f_n \land Id_{S^1})$,
for all $n\geq 0$. A symmetric spectrum is ring spectrum if it is also endowed with 
$S_n\times S_m$--equivariant multiplication maps
$$ M_{n,m}: X_n \land X_m \to X_{n+m} $$
and unit maps $\iota_0: S^0 \to X_0$ and $\iota_1: S^1 \to X_1$. They
must satisfy the associativity laws consisting of commutative squares
$$ 
\xymatrix{ X_n \land X_m \land X_r \ar[d]_{M_{n,m}\land Id} \ar[r]^{Id \land M_{m,r}} & X_n \land X_{m+r}
\ar[d]^{M_{n,m+r}} \\
X_{n+m} \land X_r \ar[r]^{M_{n+m,r}} & X_{n+m+r} } 
$$
the unit relations
$$ M_{n,0}\circ (Id\land \iota_0)= Id: X_n \simeq X_n \land S^0 \to X_n \land X_0 \to X_n $$
and similarly $M_{0,n}\circ (\iota_0\land Id)=Id$ for all $n\geq 0$, as well as
$\chi_{n,1}\circ (M_{n,1}\circ (Id\land\iota_1)) = (M_{1,n}\circ (\iota_1\land Id))\circ \tau$ with
$\tau: X_n\land S^1 \to S^1 \land X_n$ and 
$\chi_{n,m}\in S_{n+m}$ the shuffle permutation moving the first $n$ elements past the last $m$.
Commutativity of a symmetric ring spectrum is expressed by the commutativity of the diagrams
$$ \xymatrix{ X_n \land X_m \ar[r] \ar[d]_{M_{n,m}} & X_m \land X_n \ar[d]^{M_{m,n}} \\
X_{n+m} \ar[r]_{\chi_{n,m}} & X_{m+n} } $$
with the twist as the first map.
For a detailed introduction to symmetric spectra we refer the reader to [Schw12].

\medskip

{\bf 4.2.  $\Gamma$--spaces.} We recall the setting of $\Gamma$--spaces used in the Segal's
construction of spectra from categorical data. The notion of $\Gamma$--spaces
and its relation to connective spectra formalises the intuition that spectra
are a natural homotopy-theoretic generalisation of abelian groups.

\smallskip

Let $\Gamma^0$ denote the category of finite pointed sets, with
objects $$\underline{n}=\{ 0,1,2, \ldots, n \}$$ and morphisms
$f \in \Gamma^0(\underline{n},\underline{m})$ given by functions
$$
 f: \{ 0,1,2, \ldots, n \} \to \{ 0,1,2, \ldots, m \}, \ \ \ \text{ with } \ f(0)=0. 
 $$
Let $\Gamma$ denote the opposite category.

\smallskip

A pointed category $\Cal{C}$ is a category with a chosen object that
is both initial and final. A pointed functor $F : \Gamma^0 \to \Cal{C}$
is called a $\Gamma$--object in $\Cal{C}$.

\smallskip

Given a pointed category $\Cal{C}$, the category $\Gamma\Cal{C}$ has objects
the pointed functors $F : \Gamma^0 \to \Cal{C}$ and morphisms the natural
transformations between these functors. 

\smallskip

$\Gamma$--spaces are objects of the category $\Gamma\Cal{C}$, in the
case where $\Cal{C}={\underline{\Delta}}_*$ is the category of pointed simplicial sets.

\smallskip

Given a $\Gamma$--space $F : \Gamma^0 \to {\underline{\Delta}}_*$, the morphisms
$f_j :  \underline{n} \to \underline{1}$ that map the $j$-th element to $1$ 
and the rest to $0$ induce, for each $n\geq 1$, a morphism
$$
 F(\underline{n})\to \prod_{j=1}^n F(\underline{1}) .
\eqno(4.1)
$$ 

\smallskip

The {\it special} $\Gamma$--spaces (or Segal $\Gamma$--spaces)
are $\Gamma$--spaces $F$ as 
above, where all the maps (4.1) are weak equivalences.
For special $\Gamma$--spaces  the weak equivalence
$F(\underline{2})\simeq F(\underline{1})\times F(\underline{1})$
induces a monoid  
$$
 \pi_0(F(\underline{1})) \times \pi_0(F(\underline{1}) ) \to \pi_0(F(\underline{2}))
\to \pi_0(F(\underline{1})). 
$$
Such a $\Gamma$--space is called {\it very special} when this monoid is an abelian group. 

\smallskip

The $\Gamma$--space $\bold{S}: \Gamma^0 \to {\underline{\Delta}}_*$ is given
by the inclusion of the category $\Gamma^0$ into ${\underline{\Delta}}_*$
mapping a finite pointed set to the corresponding discrete pointed 
simplicial set. As shown in [Se74] (Barratt--Priddy--Quillen theorem), 
the associated spectrum is the sphere spectrum, which we also denote by $\bold{S}$. 

\smallskip

The category $\Gamma^0$ of finite pointed sets has a smash product
functor $\land: \Gamma^0 \times \Gamma^0 \to \Gamma^0$, with 
$(\underline{n}, \underline{m})\mapsto \underline{n}\land \underline{m}$,
which extends to a smash product of arbitrary pointed (simplicial) sets. 

\smallskip

The smash product of $\Gamma$-spaces constructed in [Ly99] 
is obtained by first associating to a pair $F,F': \Gamma^0 \to {\underline{\Delta}}_*$
of $\Gamma$--spaces a bi--$\Gamma$--space 
$F\tilde\land F': \Gamma^0\times \Gamma^0\to  {\underline{\Delta}}_*$
$$ 
(F \tilde\land F') (\underline{n},\underline{m}) = F(\underline{n}) \land F'(\underline{m}) 
$$
and then defining
$$ 
(F\land F')(\underline{n}) = {\roman{colim}}_{ \underline{k}\land \underline{\ell} \to \underline{n}}
(F \tilde\land F') (\underline{k},\underline{\ell}), 
$$
where $\underline{k}\land \underline{\ell}$ is the smash product $\land: \Gamma^0\times \Gamma^0 \to \Gamma^0$. It is shown in [Ly99] that, up to natural isomorphism, this smash product is 
associative and commutative and with unit given by the $\Gamma$--space $\bold{S}$, and that
the category of $\Gamma$--spaces is symmetric monoidal with respect to this product. 

\medskip

{\bf 4.3. From $\Gamma$--spaces to connective spectra}.
The construction of [Se74], and more generally [BousFr78], 
assigns a connective spectrum to a $\Gamma$--space 
in such a way as to obtain an equivalence between the homotopy category of
$\Gamma$--spaces and the homotopy category of connective spectra. 
The construction of spectra from $\Gamma$--spaces can be performed in the modern
setting of symmetric spectra, rather than in the original Bousfield--Friedlander 
formulation of [BousFr78], see Chapter~I, Section~7.4 of [Schw12].

\smallskip

If $X$ is a simplicial set, one denotes by $X_*$ the pointed simplicial set obtained
by adding a disjoint base point. Given a $\Gamma$-space 
$F: \Gamma^0 \to {\underline{\Delta}}_*$ and a pointed simplicial set $X$,
one obtains a new $\Gamma$-space  $X \land F$, which maps $\underline{n}\in \Gamma^0$
to $X\land F(\underline{n})$ in ${\underline{\Delta}}_*$. 

\smallskip

Recall that, given a functor $F: \Cal{C}^{op}\times \Cal{C} \to \Cal{D}$, the coend
$\int^{C\in \Cal{C}} F(C,C)$ is the initial cowedge, where a cowedge to an object $X$ in $\Cal{C}$ is a
family of morphisms $h_A: A \to X$, for each $A\in \Cal{C}$, such that, for all morphisms
$f: A\to B$ in $\Cal{C}$ the following diagrams commute
$$ 
\xymatrix{  F(B,A) \ar[r]^{F(f,A)} \ar[d]^{F(B,f)} & F(A,A) \ar[d]^{h_A} \\
F(B,B) \ar[r]^{h_B} & X.
} $$

\smallskip

The key step in the construction of a connective spectrum associated to a $\Gamma$--space
consists of extending a $\Gamma$--space $F$
to an endofunctor of the category of pointed simplicial sets. This endofunctor is defined in
[BousFr78] (see also [Schw99]) as the functor (still denoted by $F$) that maps a pointed
simpliciat set $K \in {\underline{\Delta}}_*$ to the coend
$$ 
F : K \mapsto \int^{\underline{n}\in \Gamma^{op}} K^n \land F(\underline{n}), 
$$
with natural assembly maps $K\land F(K')\to F(K\land K')$. 

\smallskip

The spectrum associated to $F$, which we denote by $F(\bold{S})$, is then given
by the sequence of pointed simplicial sets $F(\bold{S})_n = F(S^n=S^1\land \cdots \land S^1)$,
with the maps $S^1 \land F(S^n) \to F(S^{n+1})$.

\smallskip

The smash product of $\Gamma$--spaces is compatible with the smash product of
spectra:  in [Ly99] it is shown that if $F,F': \Gamma^0\to {\underline{\Delta}}_*$ are
$\Gamma$--spaces with $F(\bold{S})$ and $F'(\bold{S})$ the corresponding spectra, 
then there is a map of spectra 
$$ 
F(\bold{S}) \land F'(\bold{S}) \to (F\land F')(\bold{S}) 
$$
which is natural in $(F,F')$, and a weak equivalence if one of the factors is cofibrant.

\smallskip

This gives rise to a notion of ring spectra (see [Schw99]) defined as the monoids
in the symmetric monoidal category of $\Gamma$--spaces with the smash product of
[Ly99] recalled above. One refers to these as $\Gamma$--rings. Namely, a
$\Gamma$--ring is a $\Gamma$--space $F$ endowed with unit and multiplication
maps $\bold{S}\to F$ and $F\wedge F \to F$ with associativity and unit properties 
(Sec VII.3 of [MacL71]). The associated connective spectrum of a commutative $\Gamma$--ring
is a commutative symmetric ring spectrum.  However, not all connective commutative 
symmetric ring  spectra come from a commutative $\Gamma$--ring, see [La09]. 
For a comparative view of the settings of $\Gamma$-rings and symmetric ring spectra, 
see the discussion in Section~2 of [Schw99]. 

\smallskip

If $G$ is an abelian group, there is an associated $\Gamma$-space $HG$ given on
objects by 
$$ HG(\underline{n})= G \otimes \Z[\underline{n}] \simeq G^n, $$
where $\Z[\underline{n}]$ is the free abelian group on the finite set $\underline{n}$.
If $f: \underline{n} \to \underline{m}$ is a morphism in $\Gamma^0$, then the
associated morphism $H(f): HG(\underline{n}) \to HG(\underline{m})$  maps
an $n$-tuple $(g_1,\ldots, g_n)$ (with $g_0=0$)  in $G^n$ to the $m$--tuple
$(\sum_{j\in f^{-1}(1)} g_j, \ldots, \sum_{j\in f^{-1}(n)} g_n)$. This Eilenberg--MacLane
$\Gamma$--space $HG$ maps to the Eilenberg--MacLane spectrum of $G$, which
we still denote by $HG$. If $R$ is a simplicial ring, then $HR$ is an $\bold{S}$--algebra with
multiplication $HR\wedge HR\to H(R\otimes R) \to HR$ and unit $\bold{S}\to H\Z\to HR$.

\medskip
{\bf 4.4. Assemblers, spectra, and the Grothendieck ring.} We pass now to a brief survey
of  the construction of a
spectrum associated to the Grothendieck ring of varieties developed in [Z17a] and [Z17c].

\smallskip

I.~Zakharevich developed in [Z17a] and [Z17b] 
a very general formalism for scissor--congruence relations.
The abstract form of scissor--congruence consists of
categorical data called {\it assemblers}, which in turn determine a homotopy--theoretic
{\it spectrum}, whose homotopy groups embody scissor--congruence relations.
This formalism is applied in [Z17c] in the framework producing an assembler
and a spectrum whose $\pi_0$ recovers the Grothendieck ring of varieties.
This is used to obtain a characterisation  of the kernel of multiplication by
the Lefschetz motive, which provides a general explanation for the 
observations of [Bor14], [Mart16] on the fact that the Lefschetz motive
is a zero divisor in the Grothendieck ring of varieties.

\smallskip

A {\it sieve} in a category $\Cal{C}$ is a full subcategory $\Cal{C}'$ that is
closed under precomposition by morphisms in $\Cal{C}$. A {\it Grothendieck
topology} on a category $\Cal{C}$ consists of the assignment of a 
collection $\Cal{J}(X)$ of sieves in the over category $\Cal{C}/X$, for each 
object $X$ in $\Cal{C}$, with the following  properties: 

\smallskip 
(i) the over category $\Cal{C}/X$ is
in the collection $\Cal{J}(X)$; 

\smallskip
(ii) the pullback of a sieve
in $\Cal{J}(X)$ under a morphism $f: Y\to X$ is a sieve in $\Cal{J}(Y)$;

\smallskip

(iii) given $\Cal{C}'\in \Cal{J}(X)$ and a sieve   $\Cal{T}$ in $\Cal{C}/X$, if for every
$f: Y \to X$ in $\Cal{C}'$ the pullback $f^*\Cal{T}$ is in $\Cal{J}(Y)$ then
$\Cal{T}$ is in $\Cal{J}(X)$. 

\smallskip

Let $\Cal{C}$ be a category with a Grothendieck topology.
A collection of morphisms $\{ f_i: X_i \to X \}_{i\in I}$ in $\Cal{C}$  
is a {\it covering family} if the full subcategory of $\Cal{C}/X$ 
that contains all the morphisms of $\Cal{C}$ that
factor through the $f_i$,
$$
\{ g: Y \to X\,|\, \exists i\in I\, \, h: Y\to X_i\,\, \text{ such that } f_i \circ h = g \}, 
$$
is in the sieve collection $\Cal{J}(X)$.

\smallskip

In a category $\Cal{C}$ with an initial object $\emptyset$ two morphisms $f: Y \to X$
and $g: W\to X$ are called {\it disjoint} if the pullback $Y\times_X W$ exists and is equal to $\emptyset$.
A collection $\{ f_i: X_i \to X \}_{i\in I}$ in $\Cal{C}$ is disjoint if $f_i$ and $f_j$ are disjoint
for all $i\neq j \in I$.

\smallskip

An {\it assembler category} $\Cal{C}$ is a small category endowed with
a Grothendieck topology, which has an initial object
$\emptyset$ (with the empty family as covering family), and where all 
morphisms are monomorphisms, with the property that any two finite 
disjoint covering families of $X$ in $\Cal{C}$ have a common refinement 
that is also a finite disjoint covering family. 

\smallskip

A morphism of assemblers is a functor $F: \Cal{C} \to \Cal{C}'$ that is continuous
for the Grothendieck topologies and preserves the initial object and the
disjointness property, that is, if two morphisms are disjoint in $\Cal{C}$ their images
are disjoint in $\Cal{C}'$.

\smallskip

For $X$ a finite set, the coproduct of assemblers $\bigvee_{x\in X} \Cal{C}_x$
is a category whose 
objects are the initial object $\emptyset$ and all the non--initial objects of 
the assemblers $\Cal{C}_x$. Morphisms of non--initial objects are induced by those 
of $\Cal{C}_x$. 

\smallskip

The abstract scissor congruences consist of pairs of an assembler $\Cal{C}$ and
a sieve $\Cal{D}$ in $\Cal{C}$. Given such a pair, one has an associated
assembler $\Cal{C} \smallsetminus \Cal{D}$ defined as the full subcategory of $\Cal{C}$
that contains all the objects that are not non--initial objects of $\Cal{D}$. The assembler
structure on $\Cal{C} \smallsetminus \Cal{D}$ is determined by taking as covering families
in $\Cal{C} \smallsetminus \Cal{D}$ those collections $\{ f_i : X_i \to X \}_{i\in I}$ with $X_i, X$
objects in $\Cal{C} \smallsetminus \Cal{D}$  that can be completed to a covering family in $\Cal{C}$,
namely such that there exists $\{ f_j: X_j \to X \}_{j\in J}$ with $X_j$ in $\Cal{D}$ such that
$\{ f_i : X_i \to X \}_{i\in I} \cup \{ f_j: X_j \to X \}_{j\in J}$ is a covering family in $\Cal{C}$.
There is a morphism of assemblers $\Cal{C} \to \Cal{C} \smallsetminus \Cal{D}$ that maps objects
of $\Cal{D}$ to $\emptyset$ and objects of $\Cal{C} \smallsetminus \Cal{D}$ to themselves and
morphisms with source in $\Cal{C} \smallsetminus \Cal{D}$ to themselves and morphisms
with source in $\Cal{D}$ to the unique morphism to the same target with source $\emptyset$. 
The data $\Cal{C}, \Cal{D}, \Cal{C} \smallsetminus \Cal{D}$ are an {\it abstract scissor congruence}.

\smallskip

The construction of spectra from assembler categories uses the 
general construction of spectra from categorical data is provided by the Segal
construction [Se74] of spectra from $\Gamma$--spaces, that we recalled
in section 4.2 above.

\smallskip

The main construction of [Za17a] associates to an assembler $\Cal{C}$ a homotopy--theoretic
spectrum, whose homotopy groups provide a family of associated topological
invariants satisfying versions of scissor congruence relations. The main steps of
the construction can be summarized as follows (see [Za17a]):

\smallskip

(1)  One associates to an assembler $\Cal{C}$ a category $\Cal{W}(\Cal{C})$ with objects
$\{ A_i \}_{i\in I}$, collections of non--initial objects $A_i$ of $\Cal{C}$ indexed by a finite set $I$,
and morphisms $f: \{ A_i \}_{i\in I} \to \{ B_j\}_{j\in J}$ given by a map of finite sets $f: I\to J$
and morphisms $f_i: A_i \to B_{f(i)}$ such that $\{ f_i: A_i \to B_j\, :\, i\in f^{-1}(j) \}$ is a disjoint
covering family for all $j\in J$. 

\smallskip

(2)  For a finite pointed set $(X,x_0)$ and an assembler $\Cal{C}$, one considers 
the assembler $X \wedge \Cal{C} := \bigvee_{x\in X\smallsetminus \{ x_0 \}} \Cal{C}$. The
assignment $X\mapsto \Cal{N}\Cal{W}(X\wedge \Cal{C})$, where $\Cal{N}$ is the nerve, is a
$\Gamma$--space in the sense of [Se74] recalled above, hence it defines a 
spectrum $K(\Cal{C})$ by
$$ 
X_n = \Cal{N} \Cal{W}(S^n \wedge \Cal{C}) 
$$
with structure maps $S^1 \wedge X_n \to X_{n+1}$ determined by the maps
$$ 
X\wedge \Cal{N}\Cal{W}(\Cal{C}) \to \Cal{N}\Cal{W} (X\wedge \Cal{C}). 
$$
\smallskip

(3) The group $K_0(\Cal{C}):=\pi_0 K(\Cal{C})$ is the free abelian group generated by 
objects of $\Cal{C}$ modulo the scissor--congruence relations $[A]=\sum_{i\in I} [A_i]$
for each finite disjoint covering family $\{ A_i \to A \}_{i\in I}$.

\smallskip

(4) Given a morphism $\varphi: \Cal{C}_1 \to \Cal{C}_2$ of assemblers, there is an
assembler $\Cal{C}_2/\varphi$ and a morphism $\iota: \Cal{C}_2 \to \Cal{C}_2/\varphi$ 
such that the diagram
$$
K(\iota )\circ K(\varphi ):\ K(\Cal{C}_1)\to K(\Cal{C}_2) \to K(\Cal{C}_2/\varphi )
$$
is a cofiber sequence.

\medskip

{\bf 4.5. Assembler for the equivariant Grothendieck ring.}
As we have seen in the previous section,
the equivariant version of the Grothendieck ring $K_0^{\hat\Z}(\Cal{V}_\K)$ 
 is generated by isomorphism classes of varieties with a ``good"
$\hat\Z$-action, where as before good means that each orbit is contained in an affine open subvariety 
of $X$ and that the action factors through some finite level $\Z/N\Z$. 
The scissor congruence relations in $K_0^{\hat\Z}(\Cal{V}_\K)$ are of the form
$[X]=[Y]+[X\smallsetminus Y]$ where $Y\hookrightarrow X$ and $X\smallsetminus Y \hookrightarrow X$
are $\hat\Z$-equivariant embeddings. The product is given by the Cartesian product 
endowed with the induced diagonal $\hat\Z$-action.

\medskip

{\bf 4.5.1. Lemma.} {\it
The category $\Cal{C}^{\hat\Z}$ with 
objects that are varieties $X$ with a good $\hat\Z$-action  
and morphisms that are equivariant locally closed embeddings,
endowed with the Grothendieck topology generated by the covering families 
$\{ Y \hookrightarrow X, X\smallsetminus Y \hookrightarrow X \}$ 
of $\hat\Z$-equivariant embeddings, is an assembler category.
The spectrum $K^{\hat\Z}(\Cal{V})$ determined by the assembler $\Cal{C}^{\hat\Z}$ has
$\pi_0$ given by the equivariant Grothendieck ring $K_0^{\hat\Z}(\Cal{V})$.
}

\medskip

{\bf Proof.} The first part of the statement follows as in  Example 1 of  sec.~1  of [Za17c]. The
empty set is the initial object. Finite disjoint covering families are $\hat\Z$-equivariant
maps $f_i: X_i \hookrightarrow X$ where $X_i=Y_i\smallsetminus Y_{i-1}$ for a chain
of $\hat\Z$-equivariant embeddings 
$\emptyset=Y_0\hookrightarrow Y_1 \hookrightarrow \cdots \hookrightarrow Y_n=X$. 
The property that any two finite disjoint covering families have a common
refinement follows since the category has pullbacks, \cite{Zak1}.
Morphisms are compositions of closed and open $\hat\Z$-equivariant
embeddings, hence they are all monomorphisms. 
For the second part, by Theorem 2.3 of [Za17a], if $K(\Cal{C})$ is the spectrum
determined by an assembler $\Cal{C}$, then $\pi_0 K(\Cal{C})$ is generated, as an
abelian group, by the objects of $\Cal{C}$ with the scissor--congruence relations
determined by disjoint covering families. In this case this means that $K(\Cal{C}^{\hat\Z})$
is generated by the pairs $(X,\alpha)$ of a variety $X$ with a good $\hat\Z$-action
$\alpha$ with relations $[X]=[Y]+[X\smallsetminus Y]$ for the covering families given by 
$\hat\Z$-equivariant embeddings $\{ Y \hookrightarrow X, X\smallsetminus Y \hookrightarrow X \}$.
The ring structure is coming from the symmetric monoidal structure on the category of assemblers,
which induces an $E_\infty$--ring structure on the spectrum $K(\Cal{C})$. In this case it induces the
product on $K_0^{\hat\Z}(\Cal{V})$ given by the Cartesian product of varieties with the 
diagonal $\hat\Z$-action (see also Theorem 1.4 of [Ca15].
The ring structure is induced by an $E_\infty$--ring spectrum structure on $K(\Cal{C})$
which is in turn induced by a symmetric monoidal structure on the category of assembler,
cf.~[Za17a]. $\blacksquare$

Following Theorem~4.25 of [Ca15], 
the ring structure on $K_0^{\hat\Z}(\Cal{V})$ can also be seen, as in the case of the
ordinary Grothendieck ring $K_0(\Cal{V})$, as induced on $\pi_0$ by an $E_\infty$--ring 
spectrum structure obtained from the fact that the cartesian product of varieties
determines a biexact symmetric monoidal structure on $\Cal{V}$, seen as an
SW--category (a subtractive Waldhausen category).

\medskip

{\bf 4.6. Lifting the Bost--Connes algebra to spectra.} We will now show how
 to lift the maps $\sigma_n$ and $\tilde\rho_n$ of the Bost--Connes system
from the level of the Grothendieck ring $K_0^{\hat\Z}(\Cal{V})$ to the level of the
spectrum $K^{\hat\Z}(\Cal{V})$.

\medskip

{\bf 4.6.1. Proposition.} {\it
The maps $\sigma_n(X,\alpha)=(X,\alpha\circ\sigma_n)$ and $\tilde\rho_n(X,\alpha)=(X\times Z_n, \Phi_n(\alpha))$ on $K_0^{\hat\Z}(\Cal{V})$
determine endofunctors of the assembler category $\Cal{C}^{\hat\Z}$. The 
endofunctors $\sigma_n$ are compatible with the monoidal structure induced by the
Cartesian product of varieties with diagonal $\hat\Z$--action.}

\medskip

{\bf Proof.} The endofunctors $\sigma_n$ of $\Cal{C}^{\hat\Z}$ map an object $(X,\alpha)$ to
$(X,\alpha\circ \sigma_n)$ and a pair of $\hat\Z$-equivariant embeddings
$$ (Y,\alpha |_Y) \hookrightarrow (X,\alpha) \hookleftarrow (X\smallsetminus Y, \alpha |_{X\smallsetminus Y}) $$
to the pair of embedding
$$ (Y,\alpha |_Y\circ \sigma_n) \hookrightarrow (X,\alpha\circ \sigma_n) \hookleftarrow (X\smallsetminus Y, \alpha |_{X\smallsetminus Y}\circ \sigma_n). $$
This determines the functor $\sigma_n$ on both objects and morphisms of $\Cal{C}^{\hat\Z}$.
The compatibility with the monoidal structure comes from the compatibility with Cartesian
products $\sigma_n(X,\alpha)\times \sigma_n(X',\alpha')=(X\times X',(\alpha\times\alpha')\circ \Delta \circ \sigma_n)=\sigma_n ((X,\alpha)\times (X',\alpha'))$. 
The group homomorphisms $\tilde\rho_n$ of
$K_0^{\hat\Z}(\V)$ are also induced by endofunctors of $\Cal{C}^{\hat\Z}$, which map
objects by $\tilde\rho_n(X,\alpha)=(X\times Z_n, \Phi_n(\alpha))$ and pairs of $\hat\Z$--equivariant embeddings
$$ (Y,\alpha |_Y) \hookrightarrow (X,\alpha) \hookleftarrow (X\smallsetminus Y, \alpha |_{X\smallsetminus Y}) $$
to pair of embedding
$$ (Y \times Z_n, \Phi_n(\alpha |_Y) ) \hookrightarrow (X\times Z_n, \Phi_n(\alpha)) \hookleftarrow ((X\smallsetminus Y)\times Z_n, \Phi_n(\alpha |_{X\smallsetminus Y})) $$
where $\Phi_n(\alpha |_Y)=\Phi_n(\alpha) |_Y$ and $\Phi_n(\alpha |_{X\smallsetminus Y})=\Phi_n(\alpha) |_{X\smallsetminus Y}$.
The functors $\tilde\rho_n$, however, are not compatible with the monoidal structure, and this
reflects the fact that they only induce group homomorphisms on $K_0^{\hat\Z}(\V)$ rather than
ring homomorphisms.
$\blacksquare$

\smallskip

One can obtain a similar argument working with subtractive Waldhausen categories 
as in [Ca15] in place of assemblers as in [Za17a].

\medskip

{\bf 4.7. The Kontsevich--Tschinkel Burnside ring.}
In a similar way, instead of working with the Grothendieck ring $K_0^{\hat\Z}(\V)$, we can
consider the refinement of the Grothendieck ring constructed in [KoTsch17]. We discuss
here briefly how to adapt the previous construction to this case.

\smallskip

In [KoTsch17] a refinement of the Grothendieck ring of varieties is
introduced, which is based on birational equivalence. More precisely,
for $\K$ a field of characteristic zero, the Burnside semiring ${\Burn}_+(\K)$
is defined as the set of equivalence classes of
smooth $\K$-varieties under the $\K$-birational equivalence relation,
with addition and multiplication are given by disjoint union and product
over $\K$ (Definition~2 of [KoTsch17]). The Burnside ring ${\Burn}(\K)$
is the Grothendieck ring generated by the semiring ${\Burn}_+(\K)$.
Equivalently, the Burnside ring ${\Burn}(\K)$ is generated by isomorphism
classes $[X]$ of smooth varieties over $\K$ with the equivalence relation
$[X]=[U]$ for $U\hookrightarrow X$ an open embedding with dense image.

\smallskip

To construct an assembler and an associated spectrum that
recovers the Burnside ring ${\Burn}(\K)$ as its zeroth
homotopy group, we proceed again as in [Za17a].

\medskip

{\bf 4.7.1. Lemma.} {\it
Let $\Cal{C}_{{\Burn}}$ be the category with non--initial objects given by the smooth $\K$-varieties $X$
and morphisms given by the open embeddings $U\hookrightarrow X$ with dense image. Consider the
Grothendieck topology which is generated by the open dense embeddings $U\hookrightarrow X$. The
category $\Cal{C}_{{\Burn}}$ is an assembler and the associated spectrum $K(\Cal{C}_{{\Burn}})$ has
$\pi_0 K(\Cal{C}_{{\Burn}})={\Burn}(\K)$.
}

\medskip
{\bf Proof.} The initial object is the empty scheme. If $X$ is
irreducible, a disjoint covering family consists of a single 
dense open set $U \hookrightarrow X$ and the common refinement
of two disjoint covering families $U_1  \hookrightarrow X$ and
$U_2 \hookrightarrow X$ is the dense open set $U_1\cap U_2 \hookrightarrow X$.
Morphisms are monomorphisms given by compositions of open dense embeddings.
This shows that the category $\Cal{C}_{{\Burn}}$ is an assembler. As an abelian
group, $\pi_0 K(\Cal{C}_{{\Burn}})$ is generated by the objects of $\Cal{C}_{{\Burn}}$
with relations $[X]=\sum_i [X_i]$ for $\{ f_i : X_i \to X \}$ a finite disjoint covering family.
In this case this means identifying $[X]=[U]$ for any dense open embedding
$U\hookrightarrow X$, which is the equivalence relation of ${\Burn}(\K)$.
$\blacksquare$

\medskip

It is shown in [KoTsch17] that the Burnside ring ${\Burn}(\K)$ with
the grading given by the transcendence degree, maps surjectively to the
associated graded object ${\roman{gr}}\, K_0(\Cal{V}_\K)$ with respect to the filtration of
$K_0(\Cal{V}_\K)$ by dimension
$$
{\Burn}(\K) \to \roman{gr}\,K_0(\Cal{V}_\K).
\eqno(4.2)
$$

\smallskip

As we did in the case of the Grothendieck ring, we can also consider an
equivariant version of the Kontsevich--Tschinkel Burnside ring ${\Burn}(\K)$
with respect to the group $\hat\Z$, see sec.~5 of [KoTsch17]. The corresponding
assembler and spectrum are obtained as a modification of the case discussed above.
The following statement can be proved  by arguments as in Lemma 4.7.1 and Lemma 4.5.1.

\medskip
{\bf 4.7.2. Lemma.}{\it
Let ${\Burn}^{\hat\Z}(\K)$ be generated by equivalence classes of
smooth $\K$-varieties with a good $\hat\Z$-action with respect to the
equivalence relation $[X]=[U]$ for $U\hookrightarrow X$ a $\hat\Z$-equivariant
dense open embedding. The category $\Cal{C}^{\hat\Z}_{{\Burn}}$ with objects 
the smooth $\K$--varieties $X$ with a good $\hat\Z$--action and 
and morphisms the $\hat\Z$--equivariant dense open embeddings $U\hookrightarrow X$
is an assembler with $\pi_0 K(\Cal{C}^{\hat\Z}_{{\Burn}})={\Burn}^{\hat\Z}(\K)$.
}

\medskip

We refer to $K(\Cal{C}^{\hat\Z}_{{\Burn}})$ as the $\hat\Z$--equivariant Burnside spectrum.

\medskip

The notion of an epimorphic assembler with a sink was introduced in Section~4 of [Za17a].
It denotes an assembler $\Cal{C}$ with a sink object $S$ such that $\Hom(X,S)\neq \emptyset$ 
for all other objects $X\in \Cal{C}$, and with the properties that all morphisms 
$f: X\to Y$ in $\Cal{C}$ with non--initial $X$ are epimorphisms with the set 
$\{ f: X\to Y \}$ a covering family, and for $X,Y\neq \emptyset$ no two morphisms
$X\to Z$ and $Y\to Z$ are disjoint. There is a group $G_\Cal{C}$ associated to
epimorphic assembler with a sink, whose  elements  are the equivalence classes of
pairs of morphisms $f_1,f_2: X\to S$ from a non--initial object to the sink,
where the equivalence $[f_1,f_2: X\to S]=[g_1,g_2: Y \to S]$ is determined by the existence
of an object $Z$ and maps $h_X: Z\to X$ and $h_Y: Z\to Y$ such that the following diagram commutes
$$
 \xymatrix{ 
& Y\ar[dl]_{g_1} \ar[dr]^{g_2} & \\
S & Z \ar[u]^{h_Y} \ar[d]_{h_X} & S \\
& X \ar[ul]^{f_1} \ar[ur]_{f_2} 
} 
$$
Composition of morphisms is given by any (equivalent) completion to a commutative
diagram of the form 
$$ \xymatrix{ & & W \ar[dl]^{h_1} \ar[dr]_{h_2} & & \\
& X \ar[dl]^{f_1} \ar[dr]_{f_2} & & Y \ar[dl]^{g_1} \ar[dr]_{g_2} & \\
S & & S & & S } $$
It is shown in Theorem~4.8 of [Za17a] that any choice of a morphism $f_X:X \to S$
from each object of $\Cal{C}$ to the sink object $S$ determines a morphism of
assemblers $\Cal{C} \to \bold{S}_G$, where $\bold{S}_G$ is the assembler with
objects $\emptyset$ and $\star$, a non-invertible morphism $\emptyset \to \star$ and
invertible morphisms $\Aut(\star)=G$, which has spectrum $K(\bold{S}_G)=\Sigma^\infty_+ BG$
(see Example~3.2 of [Za17a]). This morphism of assemblers $\Cal{C} \to \bold{S}_G$ 
induces an equivalence on $K$--theory. 

\medskip

{\bf 4.7.3. Lemma.} {\it The assembler $\Cal{C}^{\hat\Z}_{{\Burn}}$ is a coproduct of
epimorphic assemblers with sinks
$$ 
\Cal{C}^{\hat\Z}_{{\Burn}} = \bigvee_{[X,\alpha]} \Cal{C}^{\hat\Z}_{{\Burn}}(X,\alpha) 
$$
where 
$$
K(\Cal{C}^{\hat\Z}_{{\Burn}}(X,\alpha))\simeq \Sigma^\infty_+ B \Aut^{\hat\Z}(\K(X,\alpha)) .
$$
Here $\Aut^{\hat\Z}(\K(X,\alpha))$ is the group of $\hat\Z$--equivariant birational automorphisms 
of $X$ with a good $\hat\Z$--action $\alpha$.}

\medskip
{\bf Proof.} For an irreducible smooth projective variety $X$ with a good $\hat\Z$--action
$\alpha: \hat\Z \times X \to X$, consider the assembler 
$\Cal{C}^{\hat\Z}_{{\Burn}}(X,\alpha)$ whose objects $(U,\alpha_U) \hookrightarrow (X,\alpha)$ are
the $\hat\Z$-equivariant dense open embeddings, with $\alpha, \alpha_U$
the compatible good actions of $\hat\Z$ on $X$ and $U$, respectively. 
Arguing as in Theorem~5.3 of [Za17a] for the non--eqiuvariant case,
we see that $\Cal{C}^{\hat\Z}_{{\Burn}}(X)$ satisfies the conditions of an
epimorphic assembler with sink. The associated group $G_{(X,\alpha)}^{\hat\Z}$ consists of
equivalence classes of pairs $f_1,f_2: (U,\alpha_U) \to (X,\alpha)$, and the $f_i$ are equivariant with
respect to these actions. The group $G_{(X,\alpha)}^{\hat\Z}$ is therefore given by 
the group $\Aut^{\hat\Z}(\K(X,\alpha))$ of $\hat\Z$-equivariant birational 
automorphisms of the variety with good $\hat\Z$ action $(X,\alpha)$. 
We then have $K(\Cal{C}^{\hat\Z}_{{\Burn}}(X,\alpha))\simeq 
\Sigma^\infty_+ B \Aut^{\hat\Z}(\K(X,\alpha))$.
Moreover, we can identify the assembler $\Cal{C}^{\hat\Z}_{{\Burn}}$ with the coproduct over
equivalence classes $[X,\alpha]$ of the assemblers $\Cal{C}^{\hat\Z}_{{\Burn}}(X,\alpha)$
as above, since the morphisms of $\Cal{C}^{\hat\Z}_{{\Burn}}$ between non--initial
objects come from morphisms of the $\Cal{C}^{\hat\Z}_{{\Burn}}(X,\alpha)$ and the
objects of $\Cal{C}^{\hat\Z}_{{\Burn}}$ consist of an initial object $\emptyset$ and
the non--initial objects of the $\Cal{C}^{\hat\Z}_{{\Burn}}(X,\alpha)$ for a choice of
representatives of the classes $[X,\alpha]$. 
$\blacksquare$

\medskip

The relation between the Kontsevich--Tschinkel Burnside ring ${\Burn}^{\hat\Z}(\K)$
and the equivariant Grothendieck ring $K_0^{\hat\Z}(\Cal{V}_\K)$ can then be formulated 
at the level of assemblers and spectra in a form similar to Theorem~5.2 of [Za17a],
using the same argument, adapted to the equivariant case. Let 
$\Cal{C}^{\hat\Z,(\ell)}_\K$ denote the full sub--assembler of the assembler 
$\Cal{C}^{\hat\Z}_\K$ of Lemma~4.5.1 above, consisting of varieties of dimension at most
$n$ with good $\hat\Z$-action. 

\medskip

{\bf 4.7.4. Proposition.}{\it Let $B_n^{\hat\Z}$ denote the set of birational isomorphism classes
of varieties of dimension $n$ with good $\hat\Z$-action, through $\hat\Z$-equivariant
birational isomorphisms. The coproduct assembler
$$  \Cal{C}^{\hat\Z}_{{\Burn}, n} := \bigvee_{[X,\alpha]\in B_n^{\hat\Z}}
 \Cal{C}^{\hat\Z}_{{\Burn}}(X,\alpha)  $$
 satisfies
 $$ K(\Cal{C}^{\hat\Z}_{{\Burn}, n}) \simeq 
 \bigvee_{[X,\alpha]\in B_n^{\hat\Z}} \Sigma^\infty_+ B \Aut^{\hat\Z}(\K(X,\alpha)) 
 \simeq \text{\rm hocofib}(K(\Cal{C}^{\hat\Z,(n-1)}_\K) \to K(\Cal{C}^{\hat\Z,(n)}_\K)).  $$}

\medskip

{\bf 4.8. Burnside spectrum and Bost--Connes endomorphisms.}
The same procedure we used to lift the Bost--Connes maps $\sigma_n$ and $\tilde\rho_n$ to
the Grothendieck ring $K_0^{\hat\Z}(\V)$ and the spectrum $K^{\hat\Z}(\V)$ can be
adapted to lift the same maps to the Kontsevich--Tschinkel Burnside ring ${\Burn}^{\hat\Z}(\K)$
and the spectrum $K(\Cal{C}^{\hat\Z}_{{\Burn}})$. 

\medskip

{\bf 4.8.1. Proposition.} {\it
The maps $\sigma_n$ and $\tilde\rho_n$ of the integral Bost--Connes algebra lift to
endofunctors of the assembler category $\Cal{C}^{\hat\Z}_{{\Burn}}$, with the
 $\sigma_n$ compatible with the monoidal structure induced by the Cartesian product.
These endofunctors induce the 
corresponding maps $\sigma_n$ and $\tilde\rho_n$ 
on the Kontsevich--Tschinkel Burnside ring ${\Burn}^{\hat\Z}(\K)$. 
}
\medskip

{\bf Proof.}  We argue as in Proposition 4.6.1.
The endofunctors $\sigma_n$ of $\Cal{C}^{\hat\Z}_{{\Burn}}$ map an object $(X,\alpha)$ to
$(X,\alpha\circ \sigma_n)$ and a $\hat\Z$-equivariant dense open embedding 
$$ 
(U,\alpha |_U) \hookrightarrow (X,\alpha) 
$$
to the $\hat\Z$--equivariant dense open embedding
$$ 
(U,\alpha |_U \circ \sigma_n) \hookrightarrow (X,\alpha\circ \sigma_n) . 
$$
This determines the functor $\sigma_n$ on both objects and morphisms of $\Cal{C}^{\hat\Z}_{{\Burn}}$.
As in Proposition 4.6.1 one sees the $\sigma_n$ are compatible with Cartesian products
and induce ring homomorphisms of ${\Burn}^{\hat\Z}(\K)$.
The $\tilde\rho_n$  map objects by $\tilde\rho_n(X,\alpha)=(X\times Z_n, \Phi_n(\alpha))$ and 
$\hat\Z$-equivariant dense open embeddings  $(U,\alpha |_U) \hookrightarrow (X,\alpha)$ by
$$  (U \times Z_n, \Phi_n(\alpha |_U) ) \hookrightarrow (X\times Z_n, \Phi_n(\alpha)) $$
with $\Phi_n(\alpha |_U)=\Phi_n(\alpha) |_U$. The $\tilde\rho_n$ are not compatible with
the monoidal structure and only induce group homomorphism on ${\Burn}^{\hat\Z}(\K)$.
$\blacksquare$

\bigskip

\centerline{\bf 5. Expectation values, motivic measures, and zeta functions}

\medskip
{\bf 5.1. The Bost--Connes expectation values.}
In the case of the original Bost--Connes system, one considers
representations $\pi$ of the Bost--Connes algebra (either the integral $\Cal{A}_\Z$ or
the rational $\Cal{A}_\Q=\Q[\Q/\Z]\rtimes \N$) on a Hilbert space $\Cal{H}=\ell^2(\N)$
and associates to the algebra and the representation a dynamical system,
namely the one--parameter group of automorphism $\sigma: \R \to \Aut\,(\A)$ of the $C^*$-algebra
generated by $\A_\Q$, seen as an algebra of bounded operators on $\H$. The time
evolution satisfies the covariance condition 
$$ 
\pi(\sigma_t(a)) = e^{it H} \pi(a) e^{-itH}, 
$$
where $H$ is an (unbounded) linear operator on $\H$, the Hamiltonian of the system.
In the Bost--Connes case the time evolution is determined by $\sigma_t(\mu_n)=n^{it} \mu_n$
and $\sigma_t(x)=x$ for $x\in \Q[\Q/\Z]$.  The Hamiltonian acts on the standard orthonormal
basis of $\ell^2(\N)$ as $H \epsilon_n= \log(n)\, \epsilon_n$ and the partition function
$Z(\beta)=\roman{Tr}(e^{-\beta H})$ is the Riemann zeta function, \cite{BC}. For any element
$a\in \A_\Q$ the expectation value with respect to the Bost--Connes dynamics is
then given by
$$
\langle a \rangle_\beta = \zeta(\beta)^{-1} \roman{Tr}( \pi(a) e^{-\beta H}) =\zeta(\beta)^{-1} \sum_{n\in \N}
\langle \epsilon_n , \pi(a) e^{-\beta H} \epsilon_n\rangle.
\eqno(5.1)
$$
We can similarly construct Bost--Connes expectation values associated to the noncommutative
ring $\K_0^{\hat\Z}(\V)$ defined in sec.~3.4.

\medskip
{\bf 5.2. The equivariant Euler characteristic.}
As we discussed above, the equivariant Euler characteristic $\chi: K_0^{\hat\Z}(\V_\Q) \to K_0^{\hat\Z}(\Q)$
induces a ring homomorphism $\K_0^{\hat\Z}(\V_\Q) \to \A_\Z$ where $\A_\Z$ is the integral
Bost--Connes algebra. After tensoring with $\Q$, one obtains a morphism of crossed product
algebras $\K_0^{\hat\Z}(\V_\Q) \otimes \Q = K_0^{\hat\Z}(\V_\Q)_\Q \rtimes \N \to
\A_\Q=\Q[\Q/\Z]\rtimes \N$, with $K_0^{\hat\Z}(\V_\Q)_\Q=K_0^{\hat\Z}(\V_\Q) \otimes \Q$. 

\medskip

{\bf 5.2.1. Proposition.} {\it
Let $\pi$ be a representation of the Bost--Connes algebra $\A_\Q$ on the Hilbert space $\H=\ell^2(\N)$
with $\pi(\mu_n) \epsilon_m =\epsilon_{nm}$ and $\pi(e(r))\epsilon_n =\zeta_r^n \epsilon_n$ for
$r\to \zeta_r$ an embedding of $\Q/\Z$ as the group of roots of unity in $\C^*$. Then $\pi$ determines
a one--parameter family of group homomorphism $\varphi_\beta: K_0^{\hat\Z}(\V_\Q) \to \C$,
with $\beta\in \R^*_+$, such that for all 
$[X,\alpha]\in K_0^{\hat\Z}(\V_\Q)$ the product $\zeta(\beta) \cdot \langle [X,\alpha] \rangle_\beta$,
with $\zeta(\beta)$ the Riemann zeta funciton, is a $\Z$-combination of values at roots of unity of
the polylogarithm function $\roman{Li}_\beta(x)$. 
}

\medskip

{\bf Proof.} For the generators $a=e(r)$ of $\Z[\Q/\Z]$ the expectation value (5.1) is a polylogarithm
function evaluated at a root of unity normalized by the Riemann zeta function,
$$ \langle e(r) \rangle_\beta =\zeta(\beta)^{-1}  \sum_{n\geq 1} \zeta_r^n \, n^{-\beta} = \frac{{\roman{Li}}_\beta(\zeta_r)}{\zeta(\beta)}, $$
where $\pi(e(r))\epsilon_n =\zeta_r^n \epsilon_n$ with $r\mapsto \zeta_r$ an embedding of 
$\Q/\Z$ as the roots of unity in $\C^*$. Given a representation $\pi$ of the Bost--Connes algebra,
we compose the equivariant Euler characteristic 
$\K_0^{\hat\Z}(\V_\Q) \otimes \Q = K_0^{\hat\Z}(\V_\Q)_\Q \rtimes \N \to \A_\Q=\Q[\Q/\Z]\rtimes \N$
with the Bost--Connes expectation value $\varphi(X,\alpha)=\langle \chi(X,\alpha) \rangle_\beta$.
$\blacksquare$

\medskip
{\bf 5.3. Expectation values of motivic measures.}
Other examples can be constructed using other motivic measures. For example, 
one can consider the mixed Hodge motivic measure $h: K_0(\V_\C) \to K_0(HS)$
with $h(X)=\sum_r (-1)^r [H^r_c(X,\Q)] \in K_0(HS)$. This is a refinement of the
Hodge--Deligne polynomial motivic measure $P: K_0(\V_\C) \to \Z[u,v,u^{-1},v^{-1}]$
with $P(X,u,v)=\sum_{p,q} \dim H^{p,q} u^p v^q$. 
In the case of complex varieties with a good $\hat\Z$-action that factors through a finite
quotient $\Z/n\Z$, the graded pieces $H^{p,q}$ are $\hat\Z$-modules. The equivariant
Hodge--Deligne polynomial is then defined as the polynomial in $\Z[\Q/\Z][u,v]$ given by
$P^{\hat\Z}((X,\alpha),u,v)=\sum_{p,q} E^{p,q}(X,\alpha) u^q v^q$ with $E^{p,q}(X,\alpha) 
=\sum_k (-1)^k H^{p,q}(H^k_c(X,\C))$, with the $\hat\Z$-module structure determined by the actoin $\alpha$.
The equivariant weight polynomial is given by $W^{\hat\Z}(X,w)=P^{\hat\Z}(X,w,w)$ while evaluation
at $w=1$ recovers the equivariant Euler characteristic.
The associated expectation values are then of the form
$$ \varphi_\beta (X,\alpha)=\sum_{p,q} \langle E^{p,q}(X,\alpha)  \rangle_\beta u^q v^q $$

\medskip

{\bf 5.4. Zeta functions and assemblers.}
Passing from the level of Grothendieck rings to assemblers, spectra, and $K$-theory,
as in [Za17a]--[Za17c], also provides possible methods for lifting the zeta functions at
the level of $K$--theory. One approach, currently being developed [Za18], directly uses
assemblers and the construction of a map of assemblers between the assemblers
underlying the Grothendieck ring (and its equivariant version as discussed above)
and an assembler of almost--finite $G$--sets, by mapping a variety $X$ to the 
almost--finite set $X(\bar K)$. Another approach to the lifting of zeta functions was
developed in [CaWoZa17], using \'etale cohomology and SW-categories.  
We will return to discussing zeta functions and the lifts of the Bost--Connes system
to assemblers and spectra in the second part of this paper (in preparation).

\smallskip

Following the approach being developed in [Za18], one can show that the
equivariant Euler characteristic 
$$ 
\chi^{\hat\Z}: K_0^{\hat\Z}(\Cal{V}) \to \Z[\Q/\Z] 
$$
discussed above in sections~3.1--3.2 lifts to a map of assembler by
considering, as in section~3.2, the morphism
$$ \chi^G : K_0^G(\Cal{V}) \to A(G) \to R(G) $$
with $A(G)$ the Burnside ring, for a finite group $G$, with the
equivariant Euler characteristics defined as in [Gu-Za17] as mapping
$\chi^G(X)=\sum_k [X_k]$ with $[X_k]$ the classes in $A(G)$ of the $k$-skeleta.
In the case of $\hat\Z$ one considers the completion $\hat A(\hat\Z)=\varprojlim A(\Z/n\Z)$ 
as discussed in section~3.2, where the complete Burnstein ring $\hat A(\hat\Z)$ is seen 
as the Grothendieck ring of almost--finite $\hat\Z$--sets [DrSi88]. 
\smallskip
According to [Za18], there is a construction of an assembler of almost--finite $G$--sets,
which we denote by $\Cal{AF}^G$. The equivariant Euler characteristic
$$
 \chi^G : K_0^G(\Cal{V}) \to A(G) 
 $$
then lifts to a morphism of assemblers 
$$
 \chi^G: \Cal{C}^G  \to  \Cal{AF}^G 
 $$
since the assignment of $X$ to the union of the $X_k$ as $G$--sets
and $G$--equivariant embeddings $Y\hookrightarrow X$ and
$X\smallsetminus Y \hookrightarrow X$ to the corresponding maps
of the skeleta as $G$--sets maps disjointness morphisms in the
assembler $\Cal{C}^G$ to disjoint morphisms in the assembler
$\Cal{AF}^G$. In particular, the equivariant Euler characteristic 
$$ 
\chi^{\hat\Z}: K_0^{\hat\Z}(\Cal{V}) \to \hat A(\hat\Z) 
$$
can be lifted to a morphism of assemblers 
$$ 
\chi^{\hat\Z}: \Cal{C}^{\hat\Z}  \to  \Cal{AF}^{\hat\Z}.
$$

\bigskip

\centerline{\bf 6. Dynamical $\F_1$-structures and the Bost--Connes algebra}

\medskip
{\bf 6.1. The spectrum as Euler characteristic.}
The point of view we adopt here is similar to [EbGu--Za17].
We consider the Grothendieck ring $K_0^{\Z}(\V_\C)$ of pairs $(X,f)$ 
of a complex quasi-projective variety $X$ with an automorphism $f: X \to X$, 
such that $f_*$ in homology is quasi-unipotent. The addition is given by
disjoint union and the product by the Cartesian product. 

\smallskip

The quasi-unipotent condition ensures that the spectrum of 
the induced action $f_*: H_*(X,\Z)\to H_*(X,\Z)$ is contained
in the set of roots of unity. We can then consider the spectrum
of $f_*$ as an Euler characteristic.

\medskip

{\bf 6.1.1. Lemma.} {\it
The spectrum of the induced map on homology determines 
a ring homomorphism
$$
\sigma: K_0^{\Z}(\V_\C) \to \Z[\Q/\Z] .
\eqno(6.1)
$$
}

{\bf Proof.}  To a pair $(X,f)$ we associate the spectrum of
the map $f_*: H_*(X,\Z)\to H_*(X,\Z)$, seen as a 
subset $\Sigma(f_*) \subset \Q/\Z$ of roots of unity
counted with integer multiplicities. Thus, we have
$\sigma(X,f)=\sum_{\lambda \in \Sigma(f_*)} m_\lambda \, \lambda$.
The spectrum of a tensor product is given by the set of products of
eigenvalues of the two matrices, hence the compatibility with the
ring structure of $K_0^{\Z}(\V_\C)$.
$\blacksquare$

Under suitable assumptions on the induced map on $H^*(X,\C)$ and its
Hodge decomposition, one can also consider other kinds of motivic measures 
associated to the spectrum $\Sigma$, for example generalizations 
of the Hodge--Deligne polynomial, see [EbGu--Za17].

\medskip
{\bf 6.2. Lifting the Bost--Connes algebra to dynamical $\F_1$-structures.}
Let $Z_n$ be a zero-dimensional variety with $\# Z_n(\C)=n$. Then, 
for a given $(X,f)\in K_0^{\Z}(\V_\C)$, the {\it Verschiebung}
pair $(X\times Z_n, \Phi_n(f))$ consists of the variety $X\times Z_n$ with
the automorphism $\Phi_n(f)(x,a_i)=(x,a_{i+1})$ for $i=1,\ldots, n-1$
and $\Phi_n(f)(x,a_n)=(f(x),a_1)$. 

\medskip

{\bf 6.2.1. Lemma.} {\it
The induced map in homology $\Phi_n(f)_* : H_*(X\times Z_n,\Z)\to H_*(X\times Z_n,\Z)$ is
the Verschiebung map}.

\smallskip
{\bf Proof.}
We have $H_k(Z_n,\Z)=\Z^n$ for $k=0$ and zero otherwise, hence we can identify 
$H_*(X\times Z_n,\Z)\simeq H_*(X,\Z)^{\oplus n}$. Then the action 
$\Phi_n(f)(x,a_i)=(x,a_{i+1})$ for $i=1,\ldots, n-1$ and $\Phi_n(f)(x,a_n)=(f(x),a_1)$
induces the action $\Phi_n(f)_* =V(f_*)$ in homology.
$\blacksquare$

\smallskip

{\bf 6.3. Proposition.} {\it
The maps $\sigma_n(X,f)=(X,f^n)$ and $\tilde\rho_n(X,f) = (X\times Z_n, \Phi_n(f))$
lift the maps $\sigma_n$ and $\tilde\rho_n$ of the integral Bost--Connes algebra
to the Grothendieck ring $K_0^{\Z}(\V_\C)$, compatibly with the spectrum Euler
characteristic $\sigma: K_0^{\Z}(\V_\C) \to \Z[\Q/\Z]$.
}

\medskip

{\bf Proof.} The argument is analogous to the case of $\hat\Z$-actions 
analyzed in the previous section. Because of the relations between
Frobenius $F_n(f)=f^n$ and Verschiebung $V_n(f)$ we have
$$ \sigma_n \circ \tilde\rho_n (X,f)=\sigma_n (X\times Z_n, \Phi_n(f)) =(X\times Z_n, \Phi_n(f)^n)
=(X\times Z_n, f\times 1)=(X,f)^{\oplus n} $$
$$ \tilde\rho_n \circ \sigma_n (X,f)= \tilde\rho_n (X,f^n) =(X\times Z_n, \Phi_n(f^n)) =
(X,f) \times (Z_n,\gamma), $$
with $\gamma=\Phi_n(1): a_i \mapsto a_{i+1}$ and $a_n\mapsto a_1$,
where as before we used the relation $V_n(F_n(a)b)=a V_n(b)$, which gives
$\Phi_n(f^n)=f \Phi_n(1)$. Under the spectrum Euler characteristic map
$\sigma: K_0^{\Z}(\V_\C) \to \Z[\Q/\Z]$ we then see that we have commutative
diagrams 
$$ 
\xymatrix{ K_0^{\Z}(\V_\C) \ar[r]^{\sigma} \ar[d]^{\sigma_n} & \Z[\Q/\Z] \ar[d]^{\sigma_n} \\
K_0^{\Z}(\V_\C) \ar[r]^{\sigma} & \Z[\Q/\Z] } 
$$
and
where $\sigma_n$ are ring homomorphism and $\tilde\rho_n$ are group 
homomorphisms. 
$\blacksquare$
\smallskip
Thus, we can consider a non-commutative version of the Grothendieck
ring $K_0^{\Z}(\V_\C)$.

\medskip

{\bf 6.4. Definition.} {\it
Let $\K^{\Z}_0(\V_\C)$ be the non-commutative ring generated by 
$K_0^{\Z}(\V_\C)$ together with generators $\tilde\mu_n$ and $\mu_n^*$
satisfying the relations (3.1) for all $n,m\in \N$, and (3.2)
for all $x=(X,f) \in K_0^{\Z}(\V_\C)$ and all $n\in \N$.
}

Consider then the algebra $K_0^{\Z}(\V_\C)_\Q=K_0^{\Z}(\V_\C)\otimes_\Z \Q$.
As in the $\hat\Z$--equivariant case analyzed in the previous section, the
maps $\sigma_n$ and $\tilde\rho_n$ induce endomorphisms $\sigma_n$ 
and $\rho_n$ of $K_0^{\Z}(\V_\C)_\Q$, which determine a non-commutative
semigroup crossed product algebra. The spectrum Euler characteristic (6.1)
extends to an algebra homomorphism to the rational Bost--Connes algebra.

\medskip

{\bf 6.5. Proposition.} {\it
The algebra $\K^{\Z}_0(\V_\C)\otimes_\Z \Q$ is isomorphic to a semigroup
crossed product algebra $K_0^{\Z}(\V_\C)_\Q \rtimes \N$ with the semigroup
action given by $x \mapsto n^{-1}\tilde\rho_n(x)$. The spectrum Euler characteristic (6.1)
extends to an algebra homomorphism $\sigma: \K^{\Z}_0(\V_\C)\otimes_\Z \Q \to \A_\Q$,
where $\A_\Q=\Q[\Q/\Z]\rtimes \N$ is the rational Bost--Connes algebra.
}

\medskip

{\bf Proof.} The algebra $\K^{\Z}_0(\V_\C)\otimes_\Z \Q$ is generated by
the elements of $K_0^{\Z}(\V_\C)_\Q$ and additional generators 
$\mu_n=n^{-1} \tilde \mu_n$ and $\mu_n^*$, which satisfy the relations
$$ \mu_n^*\mu_n =1, \ \ \ \ \mu_{nm}=\mu_n\mu_m, \ \ \ \ \mu_{nm}^*=\mu_n^*\mu_m^*, \ \  \forall n,m\in\N, 
\ \ \ \ \mu_n \mu_m^* = \mu_m^*\mu_n \text{ if } (n,m)=1 $$
$$ \mu_n x \mu_n^* = \rho_n(x) \ \ \text{ with } \ \  \rho_n(x) =\frac{1}{n} \tilde\rho_n (x), $$
with $\sigma_n \rho_n(x) =x$, for all $x=(X,f) \in K_0^{\Z}(\V_\C)$. 
The semigroup action in the crossed product algebra $K_0^{\Z}(\V_\C)_\Q \rtimes \N$
is given by $x \mapsto \rho_n(x) =\mu_n x \mu_n^*$, hence one obtains an identification
of these two algebras. The morphism $\sigma: \K^{\Z}_0(\V_\C)\otimes_\Z \Q \to \A_\Q$
is the map given by the spectrum Euler characteristic on elements of $K_0^{\Z}(\V_\C)$,
extended to $\Q$-coefficients, and it maps  $\chi(\mu_n)=\mu_n$ and $\chi(\mu_n^*)=\mu_n^*$.
By Proposition 6.3, it determines a homomorphism of crossed--product
algebras. 
$\blacksquare$

\medskip
{\bf 6.6. Assembler and Bost--Connes endofuntors.}
First we consider an assembler category $\Cal{C}^\Z_\C$ associated to the ring $K^{\Z}_0(\V_\C)$
and the associated $\Gamma$--space and spectrum $K^\Z(\V_\C)$ with
$\pi_0 K^\Z(\V_\C)=K^{\Z}_0(\V_\C)$, then we show that the maps $\sigma_n$
and $\tilde\rho_n$ define endofunctors of the assembler $\Cal{C}^\Z_\C$, hence they
induce maps of spectra and induced map of the homotopy groups that recover
the Bost--Connes map on $\pi_0$.

\smallskip

{\bf 6.6.1. Proposition} {\it
Let $\Cal{C}^\Z_\C$ be the following category. Its objects the pairs $(X,f)$
of a complex quasi--projective variety $X$ with an automorphism $f: X \to X$, 
such that the induced map $f_*$ in homology is quasi--unipotent. Its morphisms
$\varphi: (Y,h)\hookrightarrow (X,f)$ are given by embeddings 
$Y\hookrightarrow X$ of components preserved by the map $f$ and  $h =f|_Y$. 
This is an assembler category, and the
associated spectrum $K^\Z(\V_\C):=K(\Cal{C}^\Z_\C)$ has $\pi_0K(\Cal{C}^\Z_\C)=K^{\Z}_0(\V_\C)$.
The maps $\sigma_n$ and $\tilde\rho_n$ on $K^{\Z}_0(\V_\C)$ lift to endofunctors
of the assembler $\Cal{C}^\Z_\C$, in which $\sigma_n$ also compatible with the monoidal structure. 
}

\medskip

{\bf Proof.} The argument is similar to the $\hat\Z$-equivariant case we discussed before.
In the category $\Cal{C}^\Z_\C$ the Grothendieck topology is generated by the covering families
$\{ (X_1,f|_{X_1})\hookrightarrow (X,f), (X_2, f|_{X_2})\hookrightarrow (X,f) \}$
with $X=X_1\sqcup X_2$ and the $X_i$ are preserves by the map $f: X\to X$. 
The empty $X$ is the initial object. The finite disjoint
covering families are given by embeddings $\varphi_i: (X_i,f|_{X_i}) \hookrightarrow (X,f)$,
where the $X_i$ are unions of components preserved by the map, $f|_{X_i}=f\circ \varphi_i$. 
Any two finite disjoint families have a common refinement since the category has 
pullbacks, \cite{Zak1} and morphisms are compositions of embeddings hence monomorphisms. 
The abelian group structure on $\pi_0 K(\Cal{C}^\Z_\C)$ is determined by the
relation $(X,f)=(X_1,f|_{X_1})+(X_1,f|_{X_2})$ for
each decomposition $X=X_1\sqcup X_2$ that is preserved by
the map $f: X \to X$. The product is determined by the symmetric monoidal
structure induced by the Cartesian product. Thus, we obtain the ring $K^{\Z}_0(\V_\C)$.
The endofunctors $\sigma_n$ map objects by $\sigma_n(X,f)=(X,f^n)$
and maps pairs of embeddings with $X=X_1\sqcup X_2$
$$\{ (X_1,f|_{X_1})\hookrightarrow (X,f) \hookleftarrow (X_2, f|_{X_2}) \}$$
to pairs of embeddings
$$\{ (X_1,f^n|_{X_1})\hookrightarrow (X,f^n) \hookleftarrow (X_2, f^n|_{X_2}) \}$$
These functors are compatible with Cartesian products, hence with the monoidal structure.
The endofunctore $\tilde\rho_n$ act on objects as $\tilde\rho_n(X,f)=(X\times Z_n,\Phi_n(f))$
and map a pair of embeddings as above to the pair
$$\{ (X_1 \times Z_n,\Phi_n(f)|_{X_1})\hookrightarrow (X \times Z_n,\Phi_n(f)) \hookleftarrow (X_2 \times Z_n, \Phi_n(f)|_{X_2}) \} ,$$
where $\Phi_n(f)|_{X_i}=\Phi_n(f|_{X_i})$. The functors
$\tilde\rho_n$ are not compatible with the monoidal structure hence they induce
group homomorphisms of $\pi_0 K(\Cal{C}^\Z_\C)$.
$\blacksquare$

Note that, unlike the $\hat\Z$-equivariant cases considered in the previous sections, 
the spectrum $K(\Cal{C}^\Z_\C)$ is not so interesting topologically, since
in the assembler we are only using decompositions into connected components. The
reason for wanting only this type of scissor-congruence relations in $K^{\Z}_0(\V_\C)$
is the spectrum Euler characteristic $\sigma: K^{\Z}_0(\V_\C) \to \Z[\Q/\Z]$, which
should map a splitting $X=X_1\sqcup X_2$ compatible with $f: X\to X$ to a
corresponding splitting $H_*(X,\Z)=H_*(X_1,\Z)\oplus H_*(X_2,\Z)$ with
quasi-unipotent maps $f_*|_{X_i}$, so that the spectrum as an element
of $\Z[\Q/\Z]$ satisfies $\sigma(X,f)=\sigma(X_1,f_1)+\sigma(X_2,f_2)$. 

\bigskip

{\bf Acknowledgment.} The second named author is partially supported by
NSF grant DMS-1707882, and by NSERC Discovery Grant RGPIN-2018-04937 
and Accelerator Supplement grant RGPAS-2018-522593.  

We thank Inna Zakharevich for several useful suggestions that were incorporated
in the final version of the article.

\newpage

\centerline{\bf References}

\medskip

[At61] M.~Atiyah. {\it Characters and cohomology of finite groups.} Publ.~Math.~IH\'ES,
t.~9 (1961), pp.~23--64.

\smallskip
[Bo11a] J.~Borger. {\it The basic geometry of Witt vectors, I: The affine case.}
J. Algebra and Number Theory,  vol.~5 (2011), no.~2, pp.~231--285.

\smallskip

[Bo11b]  J.~Borger. {\it The basic geometry of Witt vectors, II: Spaces.}
Math.~Ann.,  vol.~351 (2011), no.~2, pp.~877--933.

\smallskip

[Bor14] L.~Borisov, {\it The class of the affine line is a zero divisor in the Grothendieck ring}.
arXiv:1412.6194.

\smallskip

[BoCo95] J.B.~Bost, A.~Connes. {\it Hecke algebras, type III factors and phase transitions with spontaneous symmetry breaking in number theory}. Selecta Math. (N.S.) 1 (1995) no.~3, pp.~411--457.

\smallskip

[BousFr78] A.~K.~Bousfield, E.~M.~Friedlander. {\it Homotopy theory of $\Gamma$--spaces, spectra, and
bisimplicial sets}. In ``Geometric applications of homotopy theory II", Lecture Notes in
Math. 658 (1978), pp.~ 80--130.

\smallskip

[Ca15] J.~A.~Campbell, {\it The $K$--theory spectrum of varieties}.  arXiv:1505.03136v2

\smallskip

[CaWoZa17] J.~Campbell, J.~Wolfson, I.~Zakharevich, {\it Derived $\ell$-adic zeta functions.}
arXiv:1703.09855

\smallskip

[CC14] A.~Connes, C.~Consani. {\it On the arithmetic of the BC-system}. J. Noncommut. Geom. vol.~8 (2014),
no.~3, pp.~873--945.

\smallskip
 [CCMar09] A.~Connes, C.~Consani, M.~Marcolli. {\it Fun with $\bold{F}_1$}.
J. Number Theory, vol.~129 (2009) no.~6, pp.~1532--1561.

\smallskip

[De02] P.~Deligne. {\it Cat\'egories tensorielles}, 
Mosc. Math. J. 2 (2002), no. 2, pp.~227--248. 

\smallskip

[DrSi88] A.~W.~M.~Dress, C.~Siebeneicher. {\it The Burnside ring of profinite 
groups and the Witt vector construction}. Advances in Math.,
vol.~70 (1988), no.~1, pp.~7--132.

\smallskip

[EbGu-Za17] W.~Ebeling, S.M.~Gusein--Zade. {\it Higher--order spectra, equivariant Hodge--Deligne polynomials, and Macdonald--type equations}. In ``Singularities and computer algebra", pp.~97--108, 
Springer, 2017.

\smallskip

[EKMM97] A.D.~Elmendorf, I.~Kriz, M.A.~Mandell, J.P.~May, with an appendix by M. Cole. {\it Rings, modules, and algebras in stable homotopy theory}.  Mathematical Surveys and Monographs, Vol.~47, American Mathematical Society, 1997.

\smallskip

[FrSh81]  J.~Franks, M.~Shub. {\it The existence of Morse--Smale diffeomorphisms.}
Topology, vol.~20 (1981), no.~4, pp.~273--290.

\smallskip

[Gr81] D.~Grayson. {\it $SK_1$  of an interesting principal ideal domain.} Journ. Pure
and Applied Algebra, vol.~20 (1981), pp.~157--163.

\smallskip

[Gu-Za17]  S.~M.~Gusein-Zade. {\it Equivariant analogues of the Euler characteristic
 and Macdonald type equations}. Russian Math.~Surveys, vol.~72, no.~1 (2017), pp.~1--32. 

\smallskip

[HSS00] M.~Hovey, B.~Shipley, J.~Smith. {\it Symmetric spectra.}
J. Amer. Math. Soc. vol.13 (2000), 149--208.

\smallskip

[JiMaOnSo16] S.~Jin, W.~Ma, K.~Ono, K.~Soundararajan. {\it The Riemann hypothesis
for period polynomials of Hecke eigenforms.} Proc.~Nat.~Ac.~Sci.~USA,
vol.~113 (2016), no.~10, pp.~2603--2608. Preprint arXiv:1601.03114

\smallskip

 [KoTsch17] M.~Kontsevich, Yu.~Tschinkel. {\it Specialisation of birational types}.
 arXiv:1708.05699

\smallskip

 [La09] T.~Lawson. {\it
Commutative $\Gamma$--rings do not model all commutative ring spectra}. 
Homology Homotopy Appl. 11 (2009) , no.~2, pp.~189--194.

\smallskip

[Le81] H.~Lenstra, Jr.  {\it Grothendieck groups of abelian group rings.}
Journ. of Pure and Applied Algebra, vol~20 (1981), pp.~173--193.

\smallskip
[LeBr17]  L.~Le Bruyn. {\it Linear recursive sequences and $\roman{Spec}\,\bold{Z}$
over $F_1$}. Communications in Algebra, vol.~45, no.~7 (2017), pp.~3150--3158.

http://dx.doi.org/10.180/00927872.2016.1236116

\smallskip

[Lo99]  E.~Looijenga, {\it Motivic measures}. 
S\'eminaire N. Bourbaki, 1999-2000, exp. no 874, pp.~267--297.

\smallskip

[Ly99]  M.~Lydakis. {\it Smash products and $\Gamma$-spaces}.
Math.~Proc.~Cam.~Phil.~Soc.~126 (1999), pp.~ 311--328. 
\smallskip

[MacL71] S.~MacLane. {\it Categories for the working mathematician}. Springer, 1971.
\smallskip

[Ma10] Yu.~Manin. {\it Cyclotomy and analytic geometry over $F_1$.}   In: Quanta of Maths. Conference in honor of Alain Connes.
Clay Math. Proceedings, vol. 11 (2010), pp.~385--408.
 arXiv:0809.2716. 28 pp.

\smallskip

[Ma16] Yu.~Manin. {\it Local zeta factors and geometries under $\roman{Spec}\,\bold{Z}$.}
 Izvestiya: Mathematics, vol. 80:4 ( 2016), pp.~751--758. DOI 10.1070/IM8392 .
arXiv:1407.4969

\smallskip

 [MaRe17] M.~Marcolli, Z.~Ren, {\it $q$--Deformations of statistical mechanical systems 
 and motives over finite fields}. p-Adic Numbers Ultrametric Anal. Appl. 9 (2017) no.~3, pp.~204--227. 
 
 \smallskip
 
[MaTa17] M.~Marcolli, G.~Tabuada, {\it Bost--Connes systems, categorification, quantum statistical mechanics, and Weil numbers}. J. Noncommut. Geom. 11 (2017) no. 1, pp.~1--49. 

\smallskip

[Mart16] N.~Martin. {\it The class of the affine line is a zero divisor in the Grothendieck 
ring: an improvement.} C. R. Math. Acad. Sci. Paris 354 (2016), no. 9, pp.~936--939.

\smallskip
[OnRoSp16] K.~Ono, L.~Rolen, F.~Sprung. {\it Zeta--polynomials for modular form
periods.} Preprint arXiv:1602.00752

\smallskip

[Schw99]  S.~Schwede, {\it Stable homotopical algebra and $\Gamma$--spaces}.
 Math. Proc. Cam. Phil. Soc. vol.126 (1999), pp.~329--356.
 
 \smallskip
 
 [Schw12] S.~Schwede, {\it Symmetric spectra.} Preprint 2012,
 {\tt https://www.math.uni-bonn.de/$\sim$schwede/SymSpec-v3.pdf}

\smallskip

[Se74] G.~Segal. {\it Categories and cohomology theories}. Topology, vol.~13 (1974), pp.~293--312. 

\smallskip

[ShSu75] M.~Shub, D.~Sullivan. {\it Homology theory and dynamical systems.} Topology,
vol.14 (1975), pp.~109--132.

\smallskip

[ToVa09] B.~ To\"en, M.~Vaqui\'e. {\it Au--dessous de $\roman{Spec}\bold{Z }$.} Journ.~ of K--theory: K-theory and its Applications to Algebra, Geometry, and Topology, 2009, 3 (03), pp.437--500. 10.1017/is008004027jkt048.
arXiv:0509684

\smallskip

[Ve73]  J.~L.~Verdier.  {\it Caract\'eristique d'Euler-Poincar\'e}. 
  Bull. Soc. Math. France, vol.~ 101 (1973),  pp.~ 441--445. 
  
  \smallskip
  
[Za17a] I.~Zakharevich.
{\it The $K$--theory of assemblers}. Adv. Math. 304 (2017), pp.~1176--1218.

\smallskip

[Za17b] I.~Zakharevich. {\it On $K_1$ of an assembler}. 
J. Pure Appl. Algebra 221 (2017), no. 7, pp.~1867--1898.

\smallskip
[Za17c]  I.~Zakharevich. {\it The annihilator of the Lefschetz motive}. 
Duke Math. J. 166 (2017), no. 11, pp.~1989--2022. 

\smallskip
[Za18]  I.~Zakharevich. {\it  Private communication.} July 2018.

\newpage
\bigskip

{\bf Max--Planck--Institute for Mathematics, 

Vivatsgasse 7, Bonn 53111, Germany.}

\smallskip

manin\@mpim-bonn.mpg.de

\medskip

{\bf California Institute of Technology, USA

Perimeter Institute for Theoretical Physics, Canada

University of Toronto, Canada}
\smallskip

matilde\@caltech.edu

\enddocument